\documentclass{article}

\usepackage{arxiv}
\usepackage{amsmath}
\usepackage[utf8]{inputenc} 
\usepackage[T1]{fontenc}    
\usepackage{hyperref}       
\usepackage{url}            
\usepackage{booktabs}       
\usepackage{amsfonts}       
\usepackage{nicefrac}       
\usepackage{microtype}      
\usepackage{graphicx}
\usepackage{float}
\usepackage{afterpage}
\usepackage{booktabs, calc}
\usepackage{multicol}
\usepackage{multirow}
\usepackage{indentfirst}
\usepackage{algorithm}
\usepackage{algpseudocode}

\DeclareMathAlphabet\mathbfcal{OMS}{cmsy}{b}{n} 
\usepackage{lineno}

\newcommand\norm[1]{\left\lVert#1\right\rVert}

\title{A compact fourth-order implicit-explicit Runge-Kutta type scheme for numerical solution
of the Kuramoto-Sivashinsky equation}

\author{
 Harish Bhatt  \\
Center for Student Success and Retention\\
Savannah State University\\
Savannah, GA 31404\\
  \texttt{bhatth@savannahstate.edu} \\
   \And
Abhinandan Chowdhury \\
  Department of Mathematics\\
  Savannah State University\\
 Savannah, GA 31404 \\
  \texttt{chowdhury@savannahstate.edu} \\
}

\begin{document}
\maketitle
\begin{abstract}
\indent This manuscript introduces a fourth-order Runge-Kutta based implicit-explicit scheme in time along with compact fourth-order finite difference scheme in space for the solution of one-dimensional Kuramoto-Sivashinsky equation with periodic and Dirichlet boundary conditions, respectively. 
The proposed scheme takes full advantage of method of line (\textup{MOL}) and partial fraction decomposition techniques, therefore it just need to solve two backward Euler-type linear systems at each time step to get the solution. Performance of the scheme is investigated by testing it on some test examples and by comparing numerical results with relevant known schemes. It is found that the proposed scheme is more accurate and reliable than existing schemes to solve Kuramoto-Sivashinsky equation.
\end{abstract}

\keywords{\quad Implicit-explicit scheme; \quad Kuramoto-Sivashinsky equation; \quad Partial fraction splitting technique;\quad Compact finite difference; 
\quad Exponential time differencing}

\section{Introduction}
The present study is concerned with the development and the implementation of a numerical scheme to solve well-known Kuramoto-Sivashinsky equation (\textup{KSE} - 
for brevity) \ref{eq:main_kse} with appropriate initial and boundary conditions. The KSE is the equation governing in the frame of the weakly nonlinear approximation for the shape of 
the free surface of thin film of viscous liquid falling down a vertical plane when the capillary forces are substantial. Kapitza in his pioneering work \cite{Kapit} 
first investigated the role of viscosity on the capillary flow in thin layers. Benny \cite{Benny} carried out the boundary-layer simplifications for the 
viscous case and Homsy \cite{Homsy} developed the long-wave weakly-nonlinear approximation whose consistent application allows one to derive the following dimensionless 
equation for evolution of the scaled film thickness $u$, after rescaling the variables 
\begin{equation}\label{eq:main_kse}
u_t + u u_x + \alpha u_{xx} + \beta u_{xxxx} = 0,
\end{equation}
where $\alpha$ and $\beta$ are nonzero constants. The Eq. \ref{eq:main_kse} are usually referred as \textup{KSE} which is of great fundamental interest just in the same way as its famous counterparts Korteweg-de Vries (\textup{KdV}) and Burgers equations are. The alternative form of \textup{KSE} was obtained in \cite{kura} while deriving it as a phase equation for the complex Ginzburg-Landau equation for the evolution of reaction fronts. The \textup{KSE} also describes flame-front instabilities \cite{siva_astro, siva_fluid, clav} as well as the dynamics of viscous-fluid films flowing along walls \cite{siva_mich, siva_slang}. For more details associated with the rich phenomenology of this flow, we refer the reader to \cite{alek_nako_paku} and \cite{pumir_mann_pom} for comprehensive 
reviews of the experimental and theoretical approaches, respectively.

The \textup{KSE} is a nonlinear evolution equation which is capable of demonstrating chaotic behavior in both time and space. It contains nonlinearity, 
fourth-order dissipative term $u_{xxxx}$ and second-order source term $u_{xx}$ (anti-dissipation). In fact, the \textup{KSE} 
represents the extreme case when the dispersion is negligible in 
comparison with dissipation while the \textup{KdV} equation \cite{KdV} demonstrates the other extreme situation. The sole nonlinearity in the \textup{KSE} is the convective
term, which is also known as ``eikonal'' nonlinearity \cite{chris_vel, hym_nicol} in the framework of deriving the eq. \ref{eq:main_kse}.   

Over the years, a large numbers of numerical studies have been devoted to the \textup{KSE}; readers are referred to the review paper \cite{hym_nicol} where a study consists of thorough investigation of different regimes is presented. Among the earlier works, two numerical studies are worth mentioning -- one by Frisch et al. \cite{frisch_she_thual} where a detailed multiple-scale analysis of the \textup{KSE} with $2\pi$-periodic boundary conditions is presented, and the other by Christov and Bekyarov \cite{chris_bek} where a new Fourier-Galerkin method with a complete orthonormal system of functions in $L^2(-\infty,\infty)$ is applied to the solitary wave problem for the \textup{KSE}. 
At about the same time (1992 onwards) attention has also been focused on developing theoretically sound space and time discretization. In \cite{akris_FDD}, Akrivis applied a finite difference scheme to the Eq. \ref{eq:main_kse} with periodic boundary condition. Akrivis also reported in \cite{akris_FEM} a consistent numerical approach to solve the \textup{KSE} by employing a finite element Galerkin method with extrapolated Crank-Nicolson scheme. In both cases, rigorous error analysis has been carried out 
in order to derive refined error bound. 

Over the last few decades, a computational approach based on orthogonal spline collocation (\textup{OSC}) method while seeking numerical solution to the \textup{KSE} has turned out to be quite popular. It was Manickam et al. \cite{AVM} who first attempted to solve the \textup{KSE} by using the orthogonal cubic spline collocation method in conjunction with a second-order splitting method. Later, different types of \textup{OSC} methods, such as quintic B-Spline collocation (QBSC) method \cite{RGA}, Septic B-spline collocation (SBSC) method \cite{ZMP} have been successfully implemented to seek numerical solution to the \textup{KSE}. Furthermore, a numerical scheme based on the B-spline functions is introduced in \cite{LMD} for solving the generalized Kuramoto-Sivashinsky equation (\textup{gKSE}) where one test example is devoted to discussing the nonlinear stability with convergence for eq. \ref{eq:main_kse} subjected to Gaussian initial condition.
 
In addition,  numerous other methods including discontinuous Galerkin method \cite{yan_wang}, implicit-explicit \textup{BDF} method \cite{akris_smyr}, radial basis function 
(\textup{RBF}) based mesh-free method \cite{uddin_haq_siraj}, etc. have been proposed to find the numerical solutions of the \textup{KSE}. 
In another study \cite{otom_bogo_dub} reported very lately, a lattice Boltzmann model for the \textup{KSE} is modified to achieve an enhanced level of accuracy 
and stability. Here model's enhanced stability enables one to use larger time increments which is more than enough to compensate the extra computational cost due 
to high lattice speeds -- a substantial improvement over the existing model. Another lattice Boltzmann model (LBM) with the Chapman-Enskog expansion has been proposed 
for the \textup{gKSE} in \cite{LHM} where numerically obtained results are found to be in good agreement with the analytical results. The \textup{gKSE} is 
also studied with the aid of Chebyshev spectral collocation methods in \cite{khater_temsah} where the resultant reduced system of ordinary differential 
equations has been solved by employing the implicit-explicit \textup{BDF} method depicted in \cite{akris_smyr}. Very recently, Wade et al. \cite{WAD}  have numerically studied  the KSE equation with an additional term representing the dispersive term, arising in turbulent gas flow over laminar liquid.

In this manuscript, numerical solution of the \textup{KSE} by using the fourth-order Runge-Kutta based implicit-explicit scheme in time along with compact fourth-order finite difference scheme 
in space is proposed. The proposed scheme takes full benefit of \textup{MOL} and partial fraction decomposition techniques, therefore it just need to solve two backward Euler-type linear systems at each time step to get the solution. In addition, for the efficient implementation of the scheme it just requires to compute two LU decompositions outside the time loop. Several numerical experiments on the \textup{KSE} were run in order to study an empirical convergence analysis the accuracy of the proposed scheme with other existing schemes. The numerical results exhibit that the proposed scheme provides better accuracy in most of the cases than the existing schemes considered in this manuscript.

The remainder of the paper is organized as follows. In section 2 compact fourth-order schemes
are described to approximate $u_x,~u_{xx}$ and $u_{xxxx}$. In the section 3 the fourth-order Runge-Kutta based implicit-explicit scheme is briefly explained. The linear truncation error and stability analysis of the proposed scheme are discussed in section 4. In section 5, numerical experiments are performed on \textup{KSE} to test the accuracy and reliability of the proposed scheme. The conclusions are presented in section 6.

\section{Fourth-order compact finite differencing schemes}
In order to approximate spatial derivatives in \textup{KSE}, we partitioned the computational domain $\Omega\times [0, T]=\lbrace(x,t)|\ a\leq x\leq b,\ 0\leq t\leq T\rbrace$ 
into uniform grids described by the set of nodes $\lbrace(x_i,t_j)\rbrace$, in which $x_i=a+(i-1)h,\ i=1,\cdots,N+1,~~h= \frac{(b-a)}{N}$, $t_j=jk,\ j=0,1,\cdots,M,$ and $k=\frac{T}{M},$ 
where $h$ and $k$ are spatial and temporal step sizes, respectively.   

There are several methods which are used to generate compact finite difference formula to approximate first, second and fourth-order spatial derivatives. 
Readers are referred to \cite{LSK} and references there in for more details on how to generate compact finite difference formula (\textup{CFDF}). 
In this study, the spatial derivatives in \textup{KSE} are approximated by utilizing the following \textup{CFDF}:

\subsection{Approximation of the first derivative with periodic boundary conditions}
If $u^\prime_i$ represents an approximation of the first derivative of $u(x)$ at $x_i$ then an approximation of first derivative may be written as:
\begin{flalign}
& u^\prime_{i-1}+4u^\prime_i+u^\prime_{i+1}=\frac{3}{h}(u_{i+1}-u_{i-1}),\ i=1,\cdots,N,&\label{eq:2.1}
\end{flalign}
The truncation error for formula \ref{eq:2.1} is $\mathcal{O}(h^4)$.

The matrix representation of the scheme \ref{eq:2.1} is given as:
 \begin{flalign}
& L_1{\bf U}^\prime=M_1{\bf U}, &\label{eq:2.33}
 \end{flalign}
 where 
 \begin{displaymath}
    L_1 =
     \begin{bmatrix}
	4 & 1 & 0 & \cdots & 1\\ 
	1 & 4 & 1 & \vdots & 0\\
	0 & \ddots & \ddots & \ddots & 0\\
	\vdots & & 1 & 4 & 1  \\
	1 & \cdots & 0 & 1 & 4 
     \end{bmatrix}_{N\times N},~~
    M_1 = \frac{3}{h}
     \begin{bmatrix}
	0 & 1 & 0 & \cdots & -1\\ 
	-1 & 0 & 1 & \vdots & 0\\
	0 & \ddots & \ddots & \ddots & 0\\
	\vdots & & -1 & 0 & 1  \\
	1 & \cdots & 0 & -1 & 0 
     \end{bmatrix}_{N\times N}
     {\bf U}=     \begin{bmatrix}
     u_1\\ u_2\\ \vdots \\ u_{N-1}\\ u_N 
     \end{bmatrix}_{N\times 1}.
     \end{displaymath}
     \indent Hence the fourth-order \textup{CFDF} with periodic boundary conditions for $u_{x}$ is given by:
     \begin{flalign}
     \qquad\qquad\qquad & {\bf U}^{\prime}=L^{-1}_1M_1{\bf U}.&\label{eq:2.333}
     \end{flalign} 
\subsection{Approximation of the second derivative with periodic boundary conditions}   
If $u^{\prime\prime}_i$ represents an approximation of the second derivative of $u(x)$ at $x_i$ then an approximation of second derivatives of $u(x)$ may be written as:
 \begin{flalign}
 & u^{\prime\prime}_{i-1}+10u^{\prime\prime}_{i}+u^{\prime\prime}_{i+1}=\frac{12}{h^2}(u_{i-1}-2u_i+u_{i+1}),\ i=1,2,\cdots, N.& \label{eq:2.4}
 \end{flalign}
 The matrix representation of the scheme \ref{eq:2.4} is given as:
 \begin{flalign}
& L_2{\bf U}^{\prime\prime}=M_2{\bf U}, &\label{eq:22.3}
 \end{flalign}
 where
\begin{displaymath}
 	L_2 = 
     \begin{bmatrix}
	10 & 1 & 0 & \cdots & 1\\ 
	1 & 10 & 1 & \vdots & 0\\
	0 & \ddots & \ddots & \ddots & 0\\
	\vdots & & 1 & 10& 1  \\
	1 & \cdots & 0 & 1 & 10
	\end{bmatrix}_{N\times N},  
    M_2 = \frac{12}{h^2}
     \begin{bmatrix}
	-2 & 1 & 0 & \cdots & 1\\ 
	1 & -2 & 1 & \vdots & 0\\
	0 & \ddots & \ddots & \ddots & 0\\
	\vdots & & 1 & -2 & 1  \\
	1 & \cdots & 0 & 1 & -2 
     \end{bmatrix}_{N\times N},        
 \end{displaymath}
 	Hence the fourth-order \textup{CFDF} with periodic boundary conditions for $u_{xx}$ is given by:
 	\begin{flalign}
  &{\bf U}^{\prime\prime}=L^{-1}_2M_2{\bf U}.&\label{eq:2.7}
 	\end{flalign}
 \subsection{Approximation of the fourth derivative with periodic boundary conditions}   
If $u^{(iv)}_i$ represents an approximation of the fourth derivative of $u(x)$ at $x_i$ then an approximation of fourth derivative of $u(x)$ may be obtained by replacing ${\bf U}$ in \ref{eq:22.3} by ${\bf U}^{\prime\prime}$, that is: 
\begin{flalign*}
	&L_2{\bf U}^{(iv)}=M_2{\bf U}^{\prime\prime},&
\end{flalign*}
where $L_2$ and $M_2$ are coefficient matrices defined above.

Hence the fourth-order \textup{CFDF} with periodic boundary conditions for $u_{xxxx}$ is given by:
\begin{flalign}
  &{\bf U}^{(iv)}=L^{-2}_2M^2_2{\bf U}.&\label{eq:2.77}
 	\end{flalign}
\subsection{Approximation of first derivative with Dirichlet boundary conditions} 	
 In this case, uniform grid $x_i=a+(i-1)h,\ i=1,\cdots,N,~~h= \frac{(b-a)}{N-1}$ is assumed. The standard compact finite difference formula for first derivatives of $u(x,t)$ at interior points is:
\begin{flalign}
& u^\prime_{i-1}+4u^\prime_i+u^\prime_{i+1}=\frac{3}{h}(u_{i+1}-u_{i-1}),\ i= 2,\cdots,N-1,&\label{eq:2.11}
\end{flalign}
where $u_i\approx u(x_i)$ and $u^\prime_i \approx\frac{du(x_i)}{dx}$. The truncation error for formula \ref{eq:2.11} is $\mathcal{O}(h^4)$.
At boundary, when $i=1$ use:
\begin{flalign}
& 4u^\prime_1+12u^\prime_2=\frac{3}{h}\left(-\frac{34}{9}u_1+2u_2+2u_3-\frac{2}{9}u_4\right),&\label{eq:2.22}
\end{flalign}
and when $i=N$, the formula is:
\begin{flalign}
& 4u^\prime_{N}+12u^\prime_{N-1}=\frac{3}{h}\left(\frac{34}{9}u_N-2u_{N-1}-2u_{N-2}+\frac{2}{9}u_{N-3}\right). &\label{eq:2.33}
\end{flalign}
The truncation errors in the formula \ref{eq:2.22} and \ref{eq:2.33} are $\mathcal{O}(h^4)$.\\
Writing \ref{eq:2.11}-\ref{eq:2.33} into matrix form as:
 \begin{flalign}
& L_1{\bf U}^\prime=M_1{\bf U}, &\label{eq:2.333d}
 \end{flalign}
 where 
 \begin{displaymath}
    L_1 =
     \begin{bmatrix}
	4 & 12 & 0 & \cdots & 0\\ 
	1 & 4 & 1 & \vdots & 0\\
	0 & \ddots & \ddots & \ddots & 0\\
	\vdots & & 1 & 4 & 1  \\
	0 & \cdots & 0 & 12 & 4 
     \end{bmatrix}_{N\times N}, \ \ \ \ \ \ \ \ \ 
    M_1 = \frac{3}{h}
     \begin{bmatrix}
	-\frac{34}{9} & 2 &2 & -\frac{2}{9} & 0 & \cdots & 0\\
	-1 & 0 & 1 & 0 & 0 & & \vdots \\
	0 & -1 & 0 & 1 & 0 & \ddots & \vdots \\
	0 & 0 & -1 & 0 & 1 & \ddots & 0  \\
	\vdots & \vdots & \ddots & \ddots & \ddots & \ddots \\
	 & & 0 & 0 & -1 & 0 & 1 \\
	 0 & 0 & 0 & \frac{2}{9} & -2 & -2 & \frac{34}{9}
     \end{bmatrix}_{N\times N}
     \end{displaymath}
Hence the fourth-order \textup{CFDF} with Dirichlet boundary conditions for $u_{x}$ is given by:
     \begin{flalign}
     & {\bf U}^{\prime}=L^{-1}_1M_1{\bf U}.&\label{eq:2.3333dd}
     \end{flalign}
\subsection{Approximation of second derivative with Dirichlet boundary conditions}
The standard compact finite difference formula for second derivatives of $u(x,t)$ at interior points is:
 \begin{flalign}
 & u^{\prime\prime}_{i-1}+10u^{\prime\prime}_{i}+u^{\prime\prime}_{i+1}=\frac{12}{h^2}(u_{i-1}-2u_i+u_{i+1}),\ i=2,3,\cdots, N-1,& \label{eq:2.44}
 \end{flalign}
where $u^{\prime\prime}_{i}\approx\frac{d^2u(x_i)}{dx^2}$. The truncation error for formula \ref{eq:2.44} is $\mathcal{O}(h^4)$. \\
At boundary, when $i=1$, apply:
\begin{flalign}
& 10u^{\prime\prime}_{1}+100u^{\prime\prime}_{2}=\frac{12}{h^2}(\frac{725}{72} u_1-\frac{190}{9}u_2+\frac{145}{12}u_3-\frac{10}{9}u_4+\frac{5}{72} u_5),& \label{eq:2.55}
 \end{flalign}
 and when $i=N$,
 \begin{flalign}
& 10u^{\prime\prime}_{N}+100u^{\prime\prime}_{N-1}=\frac{12}{h^2}(\frac{725}{72} u_N-\frac{190}{9}u_{N-1+}\frac{145}{12}u_{N-2}-\frac{10}{9}u_{N-3}+\frac{5}{72} u_{N-4}).& \label{eq:2.66}
 \end{flalign}
The truncation error in both of the formulae is also $\mathcal{O}(h^4)$.
 
Writing \ref{eq:2.44}-\ref{eq:2.66} in matrix form yields:
 \begin{flalign}
&L_2{\bf U}^{\prime\prime}=M_2{\bf U},&\label{eq:2.666}
 \end{flalign}
 where 
\begin{displaymath}
 	L_2 =\begin{bmatrix}
	10 & 100 & 0 & \cdots & 0\\ 
	1 & 10 & 1 & \vdots & 0\\
	0 & \ddots & \ddots & \ddots & 0\\
	\vdots & & 1 & 10 & 1  \\
	0 & \cdots & 0 & 100 & 10 
     \end{bmatrix}_{N\times N},       
    M_2 = \frac{12}{h^2}
    \begin{bmatrix}
	\frac{725}{72} & -\frac{190}{9} & \frac{145}{12} & -\frac{10}{9} & \frac{5}{72}& \cdots & 0\\
	1 & -2 & 1 & 0 & 0 & & \vdots \\
	0 & 1 & -2 & 1 & 0 & \ddots & 0 \\
	& \vdots & \vdots & \vdots &  & \\
	\vdots &  & 0 & 1 & -2 & 1 & 0  \\
	 & & 0 & 0 & 1 & -2 & 1 \\
	 0 &\cdots & \frac{5}{72}& -\frac{10}{9} & \frac{145}{12} & -\frac{190}{9} &\frac{725}{72}
     \end{bmatrix}_{N\times N},        
     \end{displaymath}	
Hence the fourth-order CFDF with Dirichlet boundary conditions for $u_{xx}$ is given by:
\begin{flalign}
     & {\bf U}^{\prime\prime}=L^{-1}_2M_2{\bf U}.&\label{eq:2.33333}
     \end{flalign}
 \subsection{Approximation of the fourth derivative with Dirichlet boundary conditions}   
If $u^{(iv)}_i$ represents an approximation of the fourth derivative of $u(x)$ at $x_i$ then an approximation of fourth derivative of $u(x)$ may be obtained by replacing $U$ in \ref{eq:2.666} by $U^{\prime\prime}$, that is: 
\begin{flalign}
	&L_2{\bf U}^{(iv)}=M_2{\bf U}^{\prime\prime},&
\end{flalign}
where $L_2$ and $M_2$ are coefficient matrices defined above.

Hence the fourth-order CFDF with Dirichlet boundary conditions for $u_{xxxx}$ is given by:
\begin{flalign}
  &{\bf U}^{(iv)}=(L^{-1}_2M_2)^2{\bf U}.&\label{eq:2.777}
 	\end{flalign}
\section{Fourth-order time stepping scheme}
This section presents brief derivation procedure of fourth-order implicit-explicit (IMEX4) scheme based on fourth-order exponential time differencing Runge-Kutta (ETDRK4-B)\cite{KRO} time integrator to solve following equation with periodic boundary conditions: 
\begin{flalign}
& \frac{\partial u}{\partial t}+\mathcal{L}u=\mathcal{F}(u,t),~~u(x,0)=g(x)&\label{eq:2.3.1}
\end{flalign} 
where $\mathcal{L}$ and $\mathcal{F}$ represent linear and nonlinear operators of KSE, respectively. 
Approximating $\mathcal{L}u = \alpha u_{xx}+ \beta u_{xxxx}$ with schemes \ref{eq:2.7} and \ref{eq:2.77} respectively and $\mathcal{F}(u,t)$ with scheme \ref{eq:2.333}, resulted the following system of ODEs:
\begin{flalign}
& \frac{\partial {\bf U}}{\partial t}+L{{\bf U}}=F({\bf U},t),& \label{eq:2.3.2}
\end{flalign} 
where $L=L^{-2}_2\big(\alpha L_2M_2+\beta M^2_2\big)$ and $F({\bf U},t)=-\frac{1}{2} L^{-1}_1M_1{\bf U}^2$.

Let $k=t_{n+1}-t_n$ be the time step size, then using a variation of constant formula, the following recurrence formula is obtained:
\begin{flalign}
\qquad\qquad &{\bf U}_{n+1}=e^{-kL}{\bf U}_n+k\int_0^1e^{-kL(1-\tau)}F({\bf U}(t_n+\tau k), t_n+\tau k)d\tau,& \label{eq:2.3.3}
\end{flalign}
where ${\bf U}_n={\bf U}(t_n)$. The expression \ref{eq:2.3.3} is an exact solution of system \ref{eq:2.3.2} and approximation of its integral term leads to various ETD schemes. This paper considers a popular ETD scheme of Runge-Kutta type, namely ETDRK4-B \cite{KRO} and presents its implicit-explicit version for its general applicability.

The ETDRK4-B \cite{KRO} scheme is given as:
 \begin{equation}
{\bf U}_{n+1}=\varphi_0(kL){\bf U}_n+k\varphi_1(kL)F_n+k\varphi_2(kL)\left(-3F_n+2F^a_n+2F^b_n-F^c_n\right)+4k\varphi_3(kL)\left( F_n-F^a_n-F^b_n+F^c_n\right), \label{2.3.6}
\end{equation}
where
\begin{flalign*}
\qquad\qquad &F_n=F({\bf U}_n,t_n), \ F^a_n=F({\bf a}_n, t_n+\frac{k}{2}),\ F^b_n=F({\bf b}_n, t_n+\frac{k}{2})\ \text{and}\ F^c_n=F({\bf c}_n, t_n+k),&\\
\qquad\qquad &{\bf a}_n=\varphi_0(kL/2){\bf U}_n+\frac{k}{2}\varphi_1(kL/2)F_n,&\\
\qquad\qquad &{\bf b}_n=\varphi_0(kL/2){\bf U}_n+\frac{k}{2}\varphi_1(kL/2)F_n+k\varphi_2(kL/2)\left( F^a_n-F_n\right),&\\
\qquad\qquad &{\bf c}_n=\varphi_0(kL){\bf U}_n+k\varphi_1(kL)F_n+2k\varphi_2(kL)\left( F^b_n-F_n\right).&\\
 \end{flalign*}
 \indent Note that the scheme \ref{2.3.6} contains matrix functions of the form:
\begin{flalign} 
\qquad\qquad & \varphi_0(kL)=e^{-kL},\ \varphi_\mu(kL)=(-kL)^{-\mu}\left(e^{-kL}-\sum_{j=0}^{\mu-1}\frac{(-kL)^j}{j!}\right),\ \mu=1,2,3.& \label{eq:2.3.4}
\end{flalign}
If an eigenvalues of $L$ close to zero, the direct computation of phi functions in \ref{eq:2.3.4} is a challenging problem in numerical analysis due 
to disastrous cancellation error in the computation. To overcome this and other numerical issues associated with it, many researchers have proposed different 
techniques to compute these phi functions. Some of the well-known techniques include (i) the rational approximation of $\varphi_\mu(kL)$ \cite{HAR,HARI} (ii) the 
Krylov subspace method \cite{HOC,TAL}, and (iii) polynomial approximation of $\varphi_\mu(kL)$ \cite{CALI,SUH}.

This study focuses on the rational approximation of $\varphi_\mu(kL)$ and introduces implicit-explicit version of ETDRK4-B schemes utilizing a partial fraction decomposition technique given in \cite{HAR}. In order to alleviate computational difficulties associated in direct computation of $\varphi_\mu(kL)$, at first a fourth-order $(2,2)$-Pad$\acute{\mathrm{e}}$ approximation to $e^{-kL}$ is utilized, which helps to avoids direct computation of the matrix exponential and higher powers of matrix inverse. In addition, another advantage we found in utilizing (2, 2)-Pad$\acute{\mathrm{e}}$ approximation is that the factors $L^{-1}$ and $L^{-3}$ cancel out in ETDRK4-B scheme.

Implementation of $(2, 2)$-Pad$\acute{\mathrm{e}}$ approximation $R_{2, 2}(kL)=(12I+6kL+k^2L^2)^{-1}(12I-6kL+k^2L^2)$ into \ref{2.3.6} to approximate matrix exponential functions yields:
\begin{equation}
{\bf U}_{n+1}=R_{2,2}(kL){\bf U}_n+P_1(kL)F_n+P_2(kL)\left(-3F_n+2F^a_n+2F^b_n-F^c_n\right)+P_3(kL)\left( F_n-F^a_n-F^b_n+F^c_n\right), \label{eq:2.3.7}
\end{equation}
where
\begin{flalign*}
\qquad\qquad\qquad & P_1(kL)=12k(12I+6kL+k^2L^2)^{-1},&\\
\qquad\qquad\qquad & P_2(kL)=k(6I+kL)(12I+6kL+k^2L^2)^{-1},&\\
\qquad\qquad\qquad & P_3(kL)=2k(4I+kL)(12I+6kL+k^2L^2)^{-1}.&
\end{flalign*}
In addition:
\begin{flalign*}
\qquad\qquad\qquad & {\bf a}_n=\tilde{R}_{2, 2}(kL){\bf U}_n+\tilde{P}_1(kL)F_n,&\\
\qquad\qquad\qquad & {\bf b}_n=\tilde{R}_{2, 2}(kL){\bf U}_n+\tilde{P}_1(kL)F_n+\tilde{P}_2(kL)\left( F^a_n-F_n\right),&\\
\qquad\qquad\qquad & {\bf c}_n=R_{2, 2}(kL){\bf U}_n+P_1(kL)F_n+2P_2(kL)\left( F^b_n-F_n\right),&
\end{flalign*}
with
\begin{flalign*}
\qquad\qquad\qquad & \tilde{R}_{2, 2}(kL)=(48I-12kL+k^2L^2)(48I+12kL+k^2L^2)^{-1},&\\
\qquad\qquad\qquad & \tilde{P}_1(kL)=24k(48I+12kL+k^2L^2)^{-1},&\\
\qquad\qquad\qquad & \tilde{P}_2(kL)=2k(12I+kL)(48I+12kL+k^2L^2)^{-1}.&
\end{flalign*}
\subsection{{Fourth-order implicit-explicit Runge-Kutta type scheme }}
\label{S:3:2}
Since the scheme \ref{eq:2.3.7} consists of high order matrix polynomials to invert, the direct implementation of it would be computationally burdensome and numerically unstable, if the matrices have high condition numbers. In addition, round off error in computing the power of the matrices can produce bad approximations \cite{CMV}. In order to handle this difficulty, ${R}_{2, 2}(kL)$ and $\tilde{R}_{2, 2}(kL)$ will not be computed directly. Instead, the problem of stably computing the inverse of matrix polynomials inherent in \ref{eq:2.3.7} is handled by utilizing a partial fraction decomposition 
technique as suggested in \cite{HAR}. This decomposition does reduce the computational complexity to just two LU decompositions over the entire time interval (provided the space step $h$, and time step $k$ are held constant). In addition, this decomposition approach is important in alleviating ill-conditioning problems because only implicit Euler type solvers are required. A description of the scheme after implementing a partial fraction decomposition technique is presented in the following algorithm and from this point the algorithm is referenced with \textup{IMEXRK4} scheme. 
%
\begin{algorithm}[H]
\caption{\textup{IMEXRK4} scheme}\label{arnoldi}
\begin{algorithmic}
\State {\bf Step 1:}
\begin{flalign*}
\text{Solve the linear system}~~~ & (kL-\tilde{c}_{1}I)R_{a}=\tilde{w_1}{\bf U}_n+k\tilde{\Omega_1}F_n.&\\
\text{Define}~~~ &{\bf a}_n={\bf U}_n+2\operatorname{Re}(R_{a}), &\\
\text{and}~~~& F^{\bf a}_n=F({\bf a}_n, t_n+\frac{k}{2}).
\end{flalign*}
\State {\bf Step 2:}
\begin{flalign*}
\text{Solve the linear system}~~~ & (kL-\tilde{c}_{1}I)R_{b}=\tilde{w_1}{\bf U}_n+k\left(\tilde{\Omega_1}-\tilde{\Omega_{2}}\right)F_n+ k\tilde{\Omega_{2}}F^{\bf a}_n,&\\
\text{Define}~~~& {\bf b}_n={\bf U}_n+2\operatorname{Re}(R_{b}),&\\
\text{and}~~~& F^{\bf b}_n=F({\bf b}_n, t_n+\frac{k}{2}).
\end{flalign*}
\State {\bf Step 3:}
\begin{flalign*}
\text{Solve the linear system}~~~ & (kL-c_{1}I)R_{c}=w_1{\bf U}_n+k\left( w_{11}-2w_{21}\right)F_n+2kw_{21}F^{\bf b}_n,&\\
\text{Define}~~~& {\bf c}_n={\bf U}_n+2\operatorname{Re}(R_{c}),&\\
\text{and}~~~& F^{\bf c}_n=F({\bf c}_n, t_n+k).
\end{flalign*}
\State {\bf Step 4:}
\begin{flalign*}
\text{Solve the linear system}~~~& (kL-c_{1}I)R_{u}=w_1{\bf U}_n+k\left( w_{11}-3w_{21}+w_{31}\right)F_n+&\\&\qquad\qquad\qquad\quad k\left(2w_{21}-w_{31}\right)\left(F^{\bf a}_n+F^{\bf b}_n\right)- k\left(w_{21}-w_{31}\right)F^{\bf c}_n.&\\
\text{Evaluate}~~~& {\bf U}_{n+1}={\bf U}_n+2\operatorname{Re}(R_{u}),&
\end{flalign*}
\end{algorithmic}
\end{algorithm}

\noindent In order to implement this IMEXRK4 scheme, poles and corresponding weights were computed for $R_{2,2}(kL),\ \lbrace P_i(kL)\rbrace_{i=1}^3,\ \tilde{R}_{2,2}(kL),\ \text{and}\ \lbrace \tilde{P}_i(kL)\rbrace_{i=1}^2$ using Maple, which are as follows:
\begin{flalign*}
& c_1 = -3.0+i1.7320508075688772935, \qquad w_1 = -6.0-i10.39230484541326376,&\\
& w_{11} = -i3.4641016151377545871,\qquad\qquad w_{21} = 0.5-i0.8660254037844386467,&\\
& w_{31} =1.0-i0.57735026918962576452,\qquad \tilde{c}_1=-6.0+i3.4641016151377545871,&\\
& \tilde{w}_1=-12.0-i20.784609690826527522,\ \quad \tilde{\Omega}_1=-i3.4641016151377545870,&\\
&\tilde{\Omega}_2=1.0-i1.7320508075688772935.& 
\end{flalign*}
\section{Linear analysis}
The linear truncation error and stability analysis of the scheme \ref{eq:2.3.7} are presented in this section.
\subsection{Truncation error analysis}
It is obvious that the overall spatial discretization is of order four because a fourth-order \textup{CFDF} is applied to linear parts of the \textup{KSE}. To analyze the overall local temporal truncation error of the scheme \ref{eq:2.3.7} for \textup{KSE}, the following linear semi-discretization system is considered:
\begin{flalign}
& \frac{\partial {\bf U}}{\partial t}=-L{\bf U}+R{\bf U},& \label{eq:trun}
\end{flalign}
where $L$ may represent matrices derived from the linear spatial discretization of second and fourth-order spatial derivatives, $R$ represents linear spatial discretization of first-order derivative term of a linear \textup{KSE} respectively, and $\bf U$ is a vector of unknowns.

Applying scheme \ref{eq:2.3.7} to Eq. \ref{eq:trun} yields:
\begin{equation}
\left.\begin{aligned}
{\bf U}_{n+1}&= \big(I+\frac{kL}{2}+\frac{k^2L^2}{12} \big)^{-1}\Big(\big(I-\frac{kL}{2}+\frac{k^2L^2}{12}+Rk\big){\bf U}_n+\frac{Rk}{2}(I+\frac{kL}{6})(-3+2{\bf a}_n+2{\bf b}_n-{\bf c}_n)+\\ &\quad  \frac{2Rk}{3}(I+\frac{kL}{4})(1-{\bf a}_n-{\bf b}_n+{\bf c}_n) \Big),&\label{eq:turn1}
\end{aligned}
\right.
\end{equation}
where
\begin{flalign*}
& {\bf a}_n=(I+\frac{1}{4}kL+\frac{1}{48}k^2L^2)^{-1}\Big(I-\frac{1}{4}kL+\frac{1}{48}k^2L^2+\frac{1}{2}Rk\Big){\bf{ U}_n},&\\
& {\bf b}_n=(I+\frac{1}{4}kL+\frac{1}{48}k^2L^2)^{-1}\Big(\big(I-\frac{1}{4}kL+\frac{1}{48}k^2L^2+\frac{1}{2}Rk\big){\bf U}_n+\frac{1}{2}Rk(I+\frac{1}{12}kL)({\bf a}_n-{\bf U}_n)\Big),&\\
& {\bf c}_n=(I+\frac{1}{2}kL+\frac{1}{12}k^2L^2)^{-1}\Big(\big(I-\frac{1}{2}kL+\frac{1}{12}k^2L^2+Rk\big){\bf U}_n+Rk(I+\frac{1}{6}kL)({\bf b}_n-{\bf U}_n)\Big).&
\end{flalign*}
Taylor expansion of Eq. \ref{eq:turn1} yields:
\begin{equation}
\left.\begin{aligned}
{\bf U}_{n+1}&=\Big(I+(R-L)k+\big(\frac{L^2}{2}-LR+\frac{R^2}{2}\big)k^2+\big(\frac{L^2R}{2}-\frac{LR^2}{2}-\frac{L^3}{6}+\frac{R^3}{6}\big)k^3+\quad\qquad\qquad\\ &\quad
\big(\frac{L^4}{24}+\frac{L^2R^2}{4}-\frac{LR^3}{6}-\frac{L^3R}{6}+\frac{R^4}{24}\big)k^4+\cdots\Big){\bf U}_n.
\end{aligned}
\right.
\end{equation}
Since the exact solution of Eq. \ref{eq:trun} is:
\begin{flalign*}
& {\bf U}(t_{n+1})=e^{(R-L)k}{\bf U}(t_n).&
\end{flalign*}
Then the local truncation error of the scheme \ref{eq:2.3.7} is:
\begin{equation}
\left.\begin{aligned}
{\bf e}_{n+1}&=\Big(I+(R-L)k+\big(\frac{L^2}{2}-RL+\frac{R^2}{2}\big)k^2+\big(\frac{RL^2}{2}-\frac{LR^2}{2}-\frac{L^3}{6}+\frac{R^3}{6}\big)k^3+\quad\qquad\qquad\\ &\quad
\big(\frac{L^4}{24}+\frac{L^2R^2}{4}-\frac{LR^3}{6}-\frac{L^3R}{6}+\frac{R^4}{24}\big)k^4+\cdots\Big){\bf U}_n-e^{(R-L)k}{\bf U}(t_n)=\mathcal{O}(k^5).
\end{aligned}
\right.
\end{equation}
Hence the scheme \ref{eq:2.3.7} is fourth-order in time discretization.
\subsection{Stability analysis}
The linear stability of the scheme \ref{eq:2.3.7} was analyzed by plotting its stability regions ( see in \cite{BEY} and references therein) for the non-linear autonomous \textup{ODE}: 
\begin{equation}
u_t=-cu+F(u), \label{eq:*}
\end{equation} where $F(u)$ is a nonlinear part and $c$ represents the approximation of the linear parts of \textup{KSE}. Let us assume that there exists a fixed point $u_0$ such that $-cu_0+F(u_0)=0.$ Linearizing about this fixed point, we thus obtain:
\begin{equation} 
u_t=-cu+\gamma u, \label{eq:**}
\end{equation}
where $u$ is perturbation of $u_0$, and $\gamma=F^{'}(u_0)$. If \ref{eq:*} represents a system of \textup{ODE}s, then $\gamma$ is a diagonal or a block diagonal matrix containing the eigenvalues of $F$. To keep the fixed point stable, we need that $\operatorname{Re}(\gamma-c)<0$, for all $\gamma$ (see\cite{COX}). This approach only provides an indication to how stable a numerical method is, since in general one cannot linearize both terms simultaneously \cite{KRO}.

In general, the parameters $c$ and $\gamma$ may both be complex-valued. The stability region of the scheme \ref{eq:2.3.7} is four dimensional and therefore difficult to represent \cite{COX}. The two-dimensional stability region is obtained, if both $c$ and $\gamma$ are purely imaginary or purely real \cite{FON}, or if $\gamma$ is complex and $c$ is fixed and real \cite{BEY}.

\begin{figure}[H]
\centerline{\includegraphics[scale=0.5]{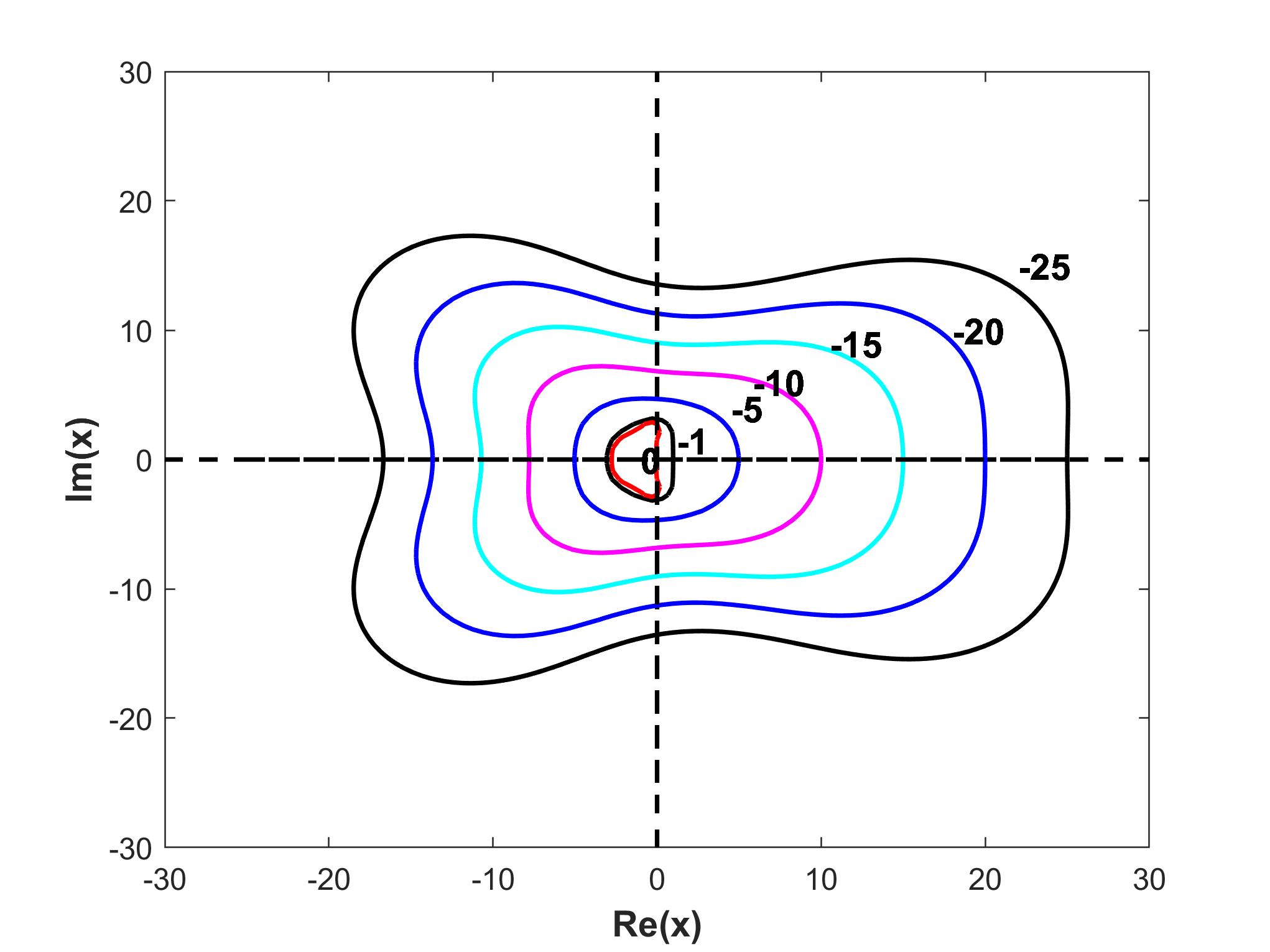}}
\caption{\footnotesize{\ Stability regions for different values of $y\in \ Re^{-}$}}
\label{fig: stability1}
\end{figure}

\begin{figure}[H]
\begin{minipage}[b]{0.5\linewidth}
\centering
\centerline{\includegraphics[scale=0.5]{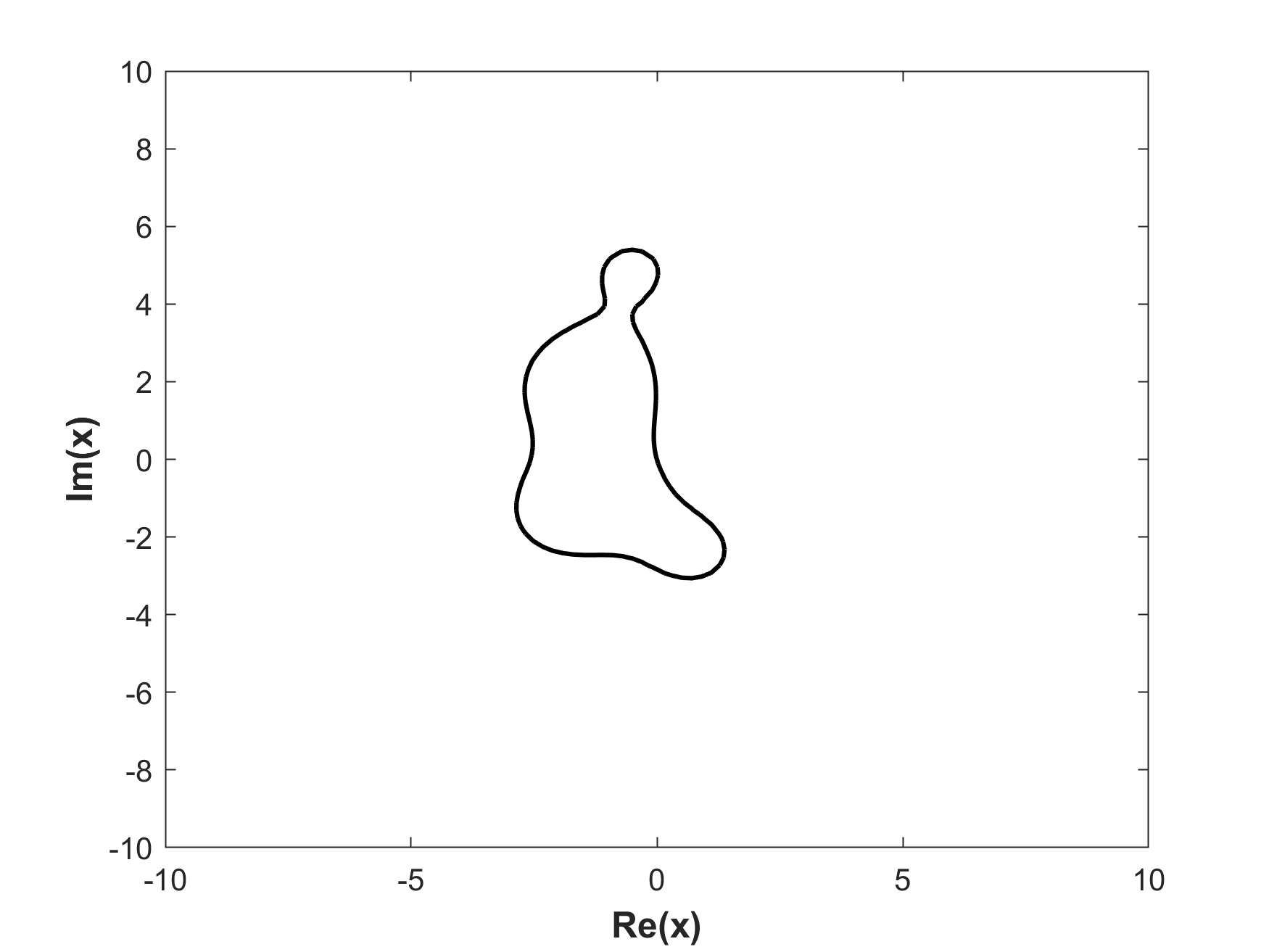}}
\centerline{\text{(a)\ $y=-5i$}}
\end{minipage}
\begin{minipage}[b]{0.5\linewidth}
\centering
\centerline{\includegraphics[scale=0.5]{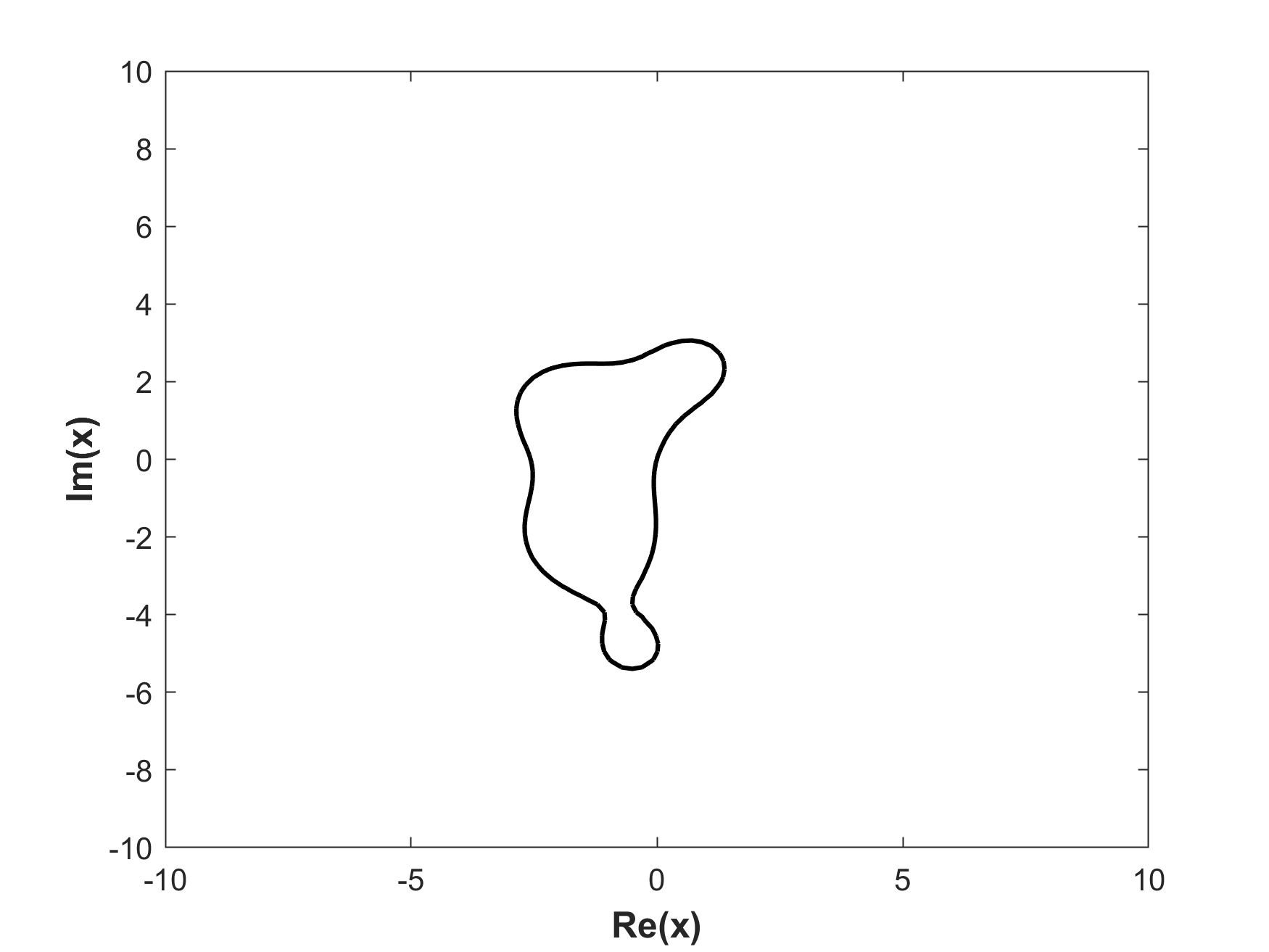}}
\centerline{\text{(b)\ $y=5i$}}
\end{minipage}

\begin{minipage}[b]{0.5\linewidth}
\centering
\centerline{\includegraphics[scale=0.5]{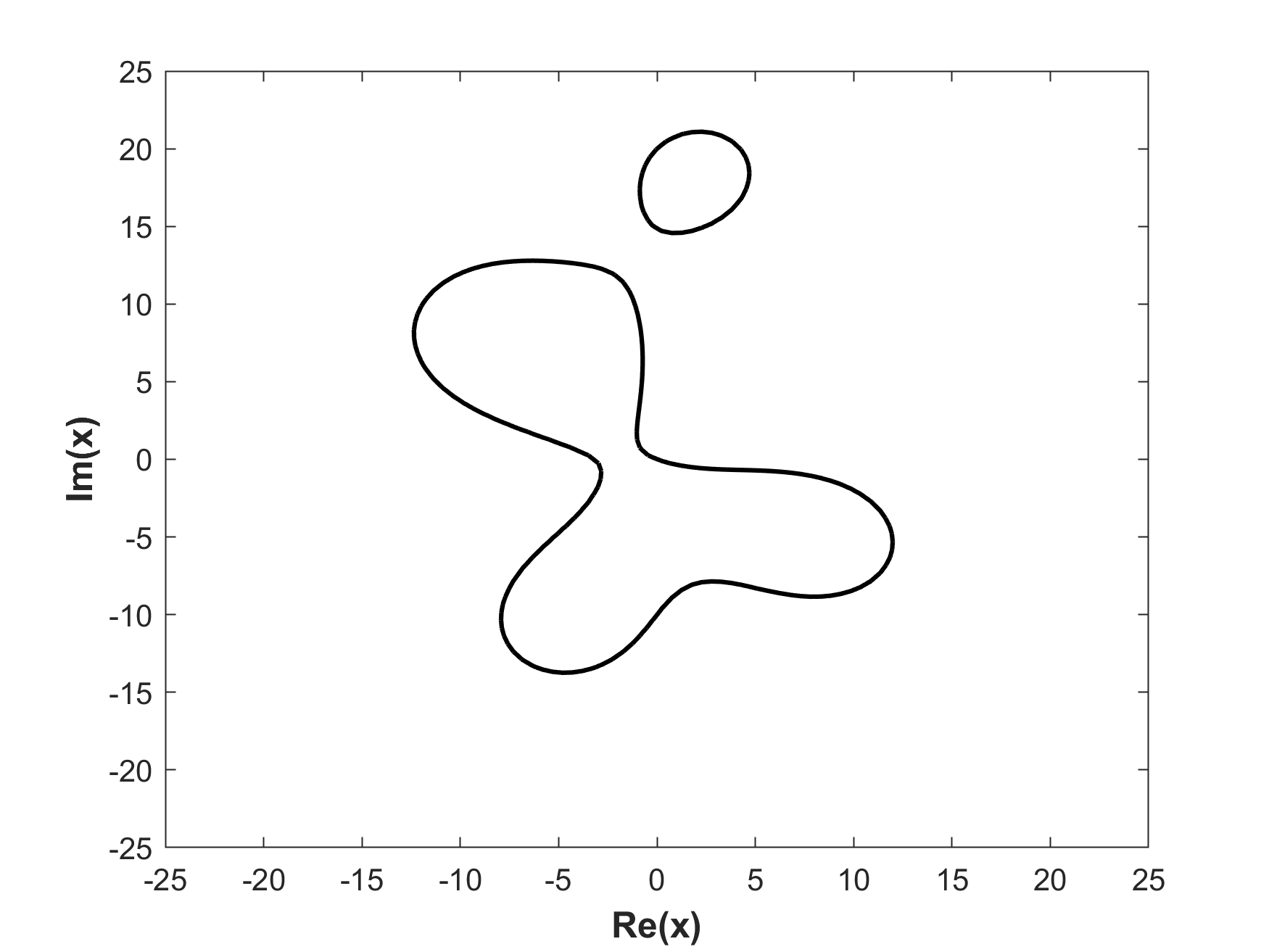}}
\centerline{\text{(c)\ $y=-20i$}}
\end{minipage}
\begin{minipage}[b]{0.5\linewidth}
\centering
\centerline{\includegraphics[scale=0.5]{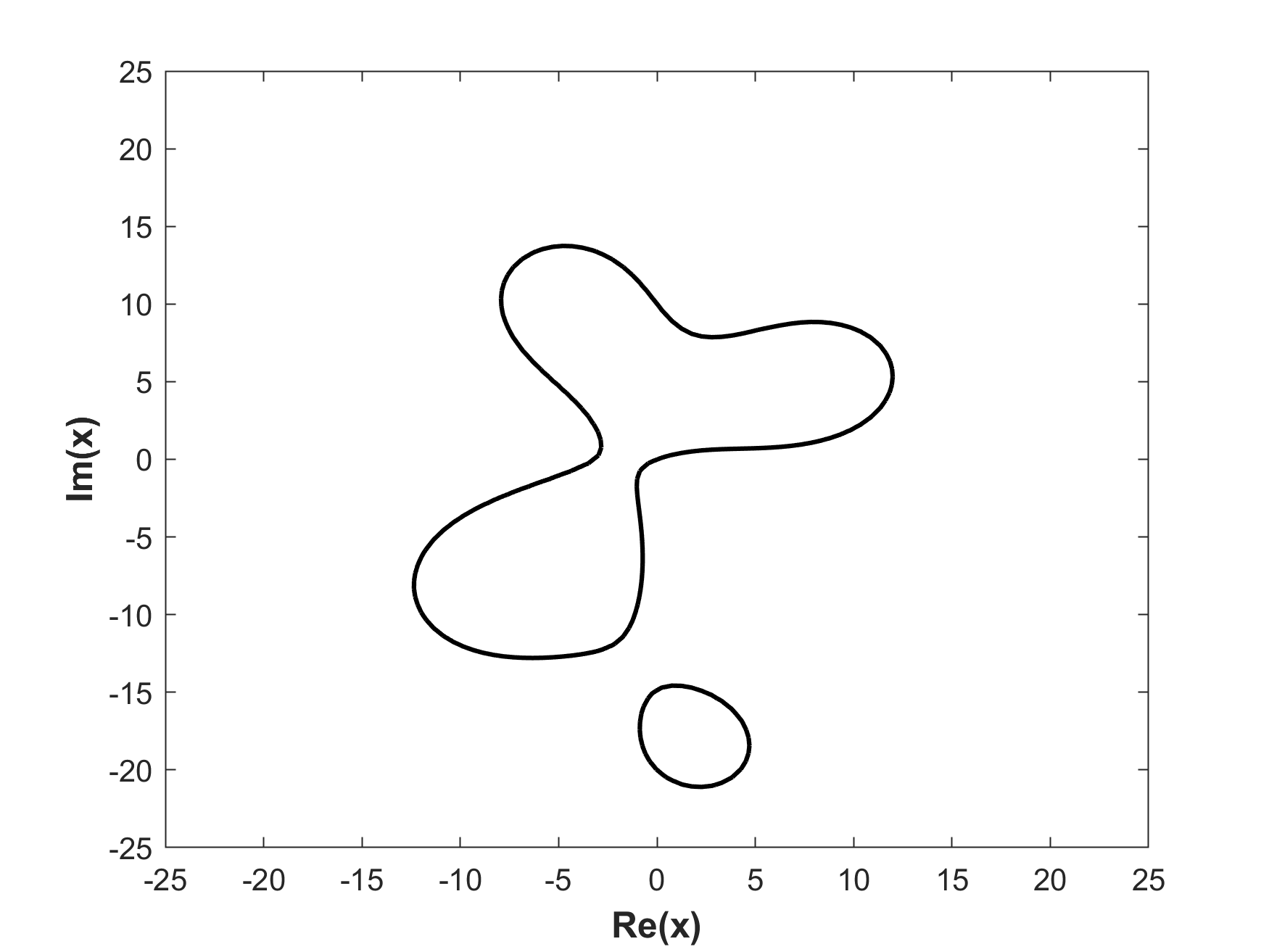}}
\centerline{\text{(d)\  $y=20i$}}
\end{minipage}
\caption{\footnotesize{\ Stability regions for different values of $y=-5i, 5i, -20i, $ and $20i$ }}
\label{fig: stability2}
\end{figure}

Utilization of the scheme \ref{eq:2.3.7} to the linearized Eq. \ref{eq:**} leads to a recurrence relation involving $u_{n}$, and $u_{n+1}.$ By letting $r=\frac{u_{n+1}}{u_n}$, $x=\gamma k$, and $y=-ck$, we come up with the following amplification factor:
\begin{flalign}
& r(x,y)= c_0+c_1x+c_2x^2+c_3x^3+c_4x^4,&
\label{eq:***}
\end{flalign} 
where
\begin{flalign*}
& c_0=1+y+\frac{1}{2}y^2+\frac{1}{6}y^3+\frac{1}{24}y^4+\frac{1}{144}y^5+O(y^6),&\\
& c_1=1+y+\frac{1}{2}y^2+\frac{1}{6}y^3+\frac{7}{192}y^4+\frac{5}{2304}y^5+O(y^6),&\\
 & c_2=\frac{1}{2}+\frac{1}{2}y+\frac{1}{4}y^2+\frac{23}{288}y^3+\frac{11}{768}y^4-\frac{37}{27648}y^5+O(y^6),&\\
 & c_3=\frac{1}{6}+\frac{1}{6}y+\frac{47}{576}y^2+\frac{77}{3456}y^3+\frac{5}{9216}y^4-\frac{515}{165888}y^5+O(y^6),&\\
 & c_4=\frac{1}{24}+\frac{1}{32}y+\frac{35}{3456}y^2-\frac{1}{13824}y^3-\frac{169}{82944}y^4-\frac{457}{331776}y^5+O(y^6).&
\end{flalign*}
The boundaries of the stability regions of the scheme \ref{eq:2.3.7} are obtained by substituting $r=e^{i\theta}, \ \theta\in [0, 2\pi]$ into the Eq. \ref{eq:***} and solving for $x$, but unfortunately we do not know the explicit expression for $|r(x,y)|=1.$  We will only be able to plot it and in this paper we have plotted the stability regions for the two cases. At first, this study focuses on the case where $\gamma$ is complex and $c$ is fixed and real. The analysis begins by selecting several real negative values of $y$ and looking for a region of stability in the complex $x$  plane where $|r(x,y)|=1$. The corresponding families of stability regions of the scheme \ref{eq:2.3.7} in the complex $x$ are plotted in Fig. \ref{fig: stability1}. According to Beylkin et al. \cite{BEY} for scheme to be  applicable, it is important that stability regions grow as $y\rightarrow -\infty$. As we can see in Fig. \ref{fig: stability1}, the stability regions for the scheme grow larger as $y\rightarrow -\infty$. These regions give an indication of the stability of the proposed scheme. 
In the second case, we assume $\gamma$ is complex and $c$ is purely imaginary and stability regions for different values of $y=-5i, 5i, -20i, $ and $20i$ are depicted in Fig.\,\ref{fig: stability2} (a)-(d).
\section{Numerical experiments and discussions}
In this section, we present the results of extensive numerical experiments carried out by implementation of the proposed scheme on four test problems 
in order to demonstrate the efficiency and accuracy of the scheme. All the numerical experiments are conducted in \textup{MATLAB} 9.3 platforms based on an Intel 
Core i52410M 2.30 GHz workstation. The accuracy of the schemes is measured in terms of maximum norm errors $\norm{\cdot}_\infty$ and global relative error (\textup{GRE}) at the final time $t=T$ 
which are defined as:
$$\norm{\cdot}_\infty=\max\limits_{1\leq i\leq N}\vert u(x_i,T)-U_i^M\vert,$$ $$ \textup{GRE}=\frac{\Sigma_i\vert u(x_i,T)-U_i^M\vert}{\Sigma_i|u(x_i,T)|},$$
where $u$, and $U$ are the exact and numerical solutions, respectively.

When the exact solution of the considered problem(s) is/are available, we compute the spatial convergence rate with: 
\begin{flalign*}
&\text{order}=\frac{\mathrm{log_{10}}(\norm{u-U_{h}}_\infty/\norm{u-U_{\frac{h}{2}}}_\infty)}{\mathrm {log_{10}}(2)}&
\end{flalign*}
and the temporal convergence rate with:
\begin{flalign*} 
&\text{order}=\frac{\mathrm{log_{10}}(\norm{u-U_{k}}_\infty/\norm{u-U_{\frac{k}{2}}}_\infty)}{\mathrm {log_{10}}(2)}&,
\end{flalign*}
where $\norm{u-U_{h}}_\infty$ and $\norm{u-U_{h/2}}_\infty$ are maximum error norms with spatial step sizes equal to $h$ and $\frac{h}{2}$, respectively. Similar description is valid for temporal convergence rate.\\
\indent On the other hand, when the exact solution of the problem(s) is/are unavailable, we utilize
\begin{flalign} 
&\text{order}=\frac{\mathrm{log_{10}}(E_k/E_{\frac{k}{2}})}{\mathrm {log_{10}}(2)},&
\end{flalign}
where $E_k=\norm{U_k-U_{2k}}_\infty$ and $E_{\frac{k}{2}}=\norm{U_{\frac{k}{2}}-U_{k}}_\infty$ are maximum norm errors at $k$ and $\frac{k}{2}$ to measure the temporal convergence rate of the scheme.

\noindent {\bf Example 1 (Benchmark Problem):}\  In this example, the Kuramoto-Sivashinsky with $\alpha=-1$ and $\beta=1$ 
over a domain $\Omega=[-50, 50]$, with the analytical solution
\begin{flalign}\label{prob1}
&u(x,t)=\mu +\frac{15\tanh^3(\nu(x-\mu t-x_0))-45\tanh(\nu(x-\mu t-x_0))}{19^{\frac{3}{2}}}
\end{flalign} is considered.

The initial and boundary conditions are extracted from the exact solution \ref{prob1}. This example is considered as a benchmark
problem in order to investigate the performance in terms of accuracy and efficiency of the proposed method. 
The parameters in \ref{prob1} are chosen as $\mu=5,~\nu=\frac{1}{2\sqrt{19}}$ and $x_0=-25$.

\begin{figure}[H]
\begin{minipage}[b]{0.5\linewidth}
\centering
\centerline{\includegraphics[width=\linewidth,height=\textheight,keepaspectratio]{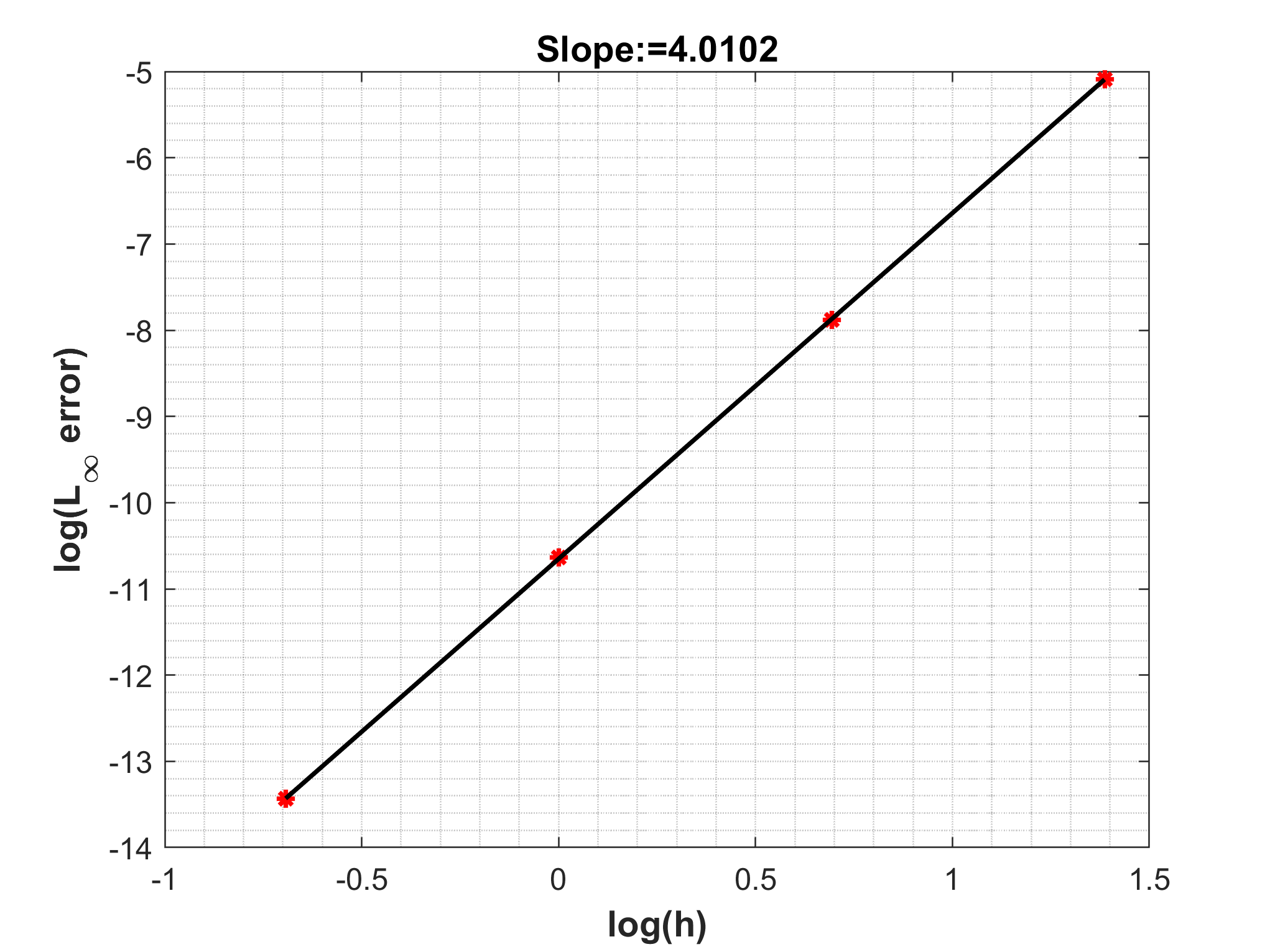}}
\centerline{\text{(a)~spatial rate of convergence}}
\end{minipage}
\begin{minipage}[b]{0.5\linewidth}
\centering
\centerline{\includegraphics[width=\linewidth,height=\textheight,keepaspectratio]{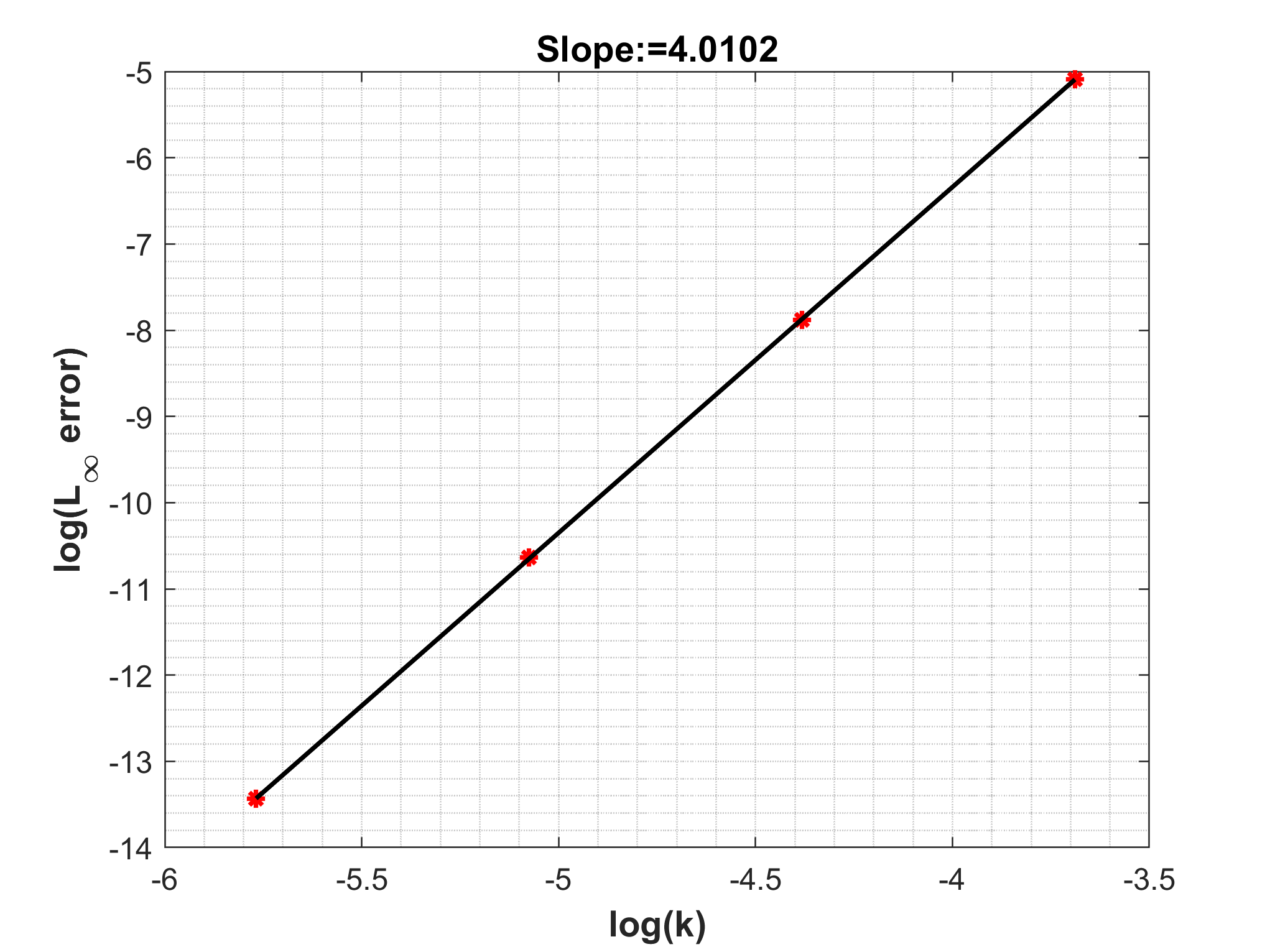}}
\centerline{\text{(b)~time rate of convergence}}
\end{minipage}
\caption{\footnotesize{Log-log plots of spatial and time rates of convergence of the proposed method.}}
\label{fig:spacetime cgce}
\end{figure}

\begin{figure}[H]
\begin{minipage}[b]{0.5\linewidth}
\centering
\centerline{\includegraphics[width=\linewidth,height=\textheight,keepaspectratio]{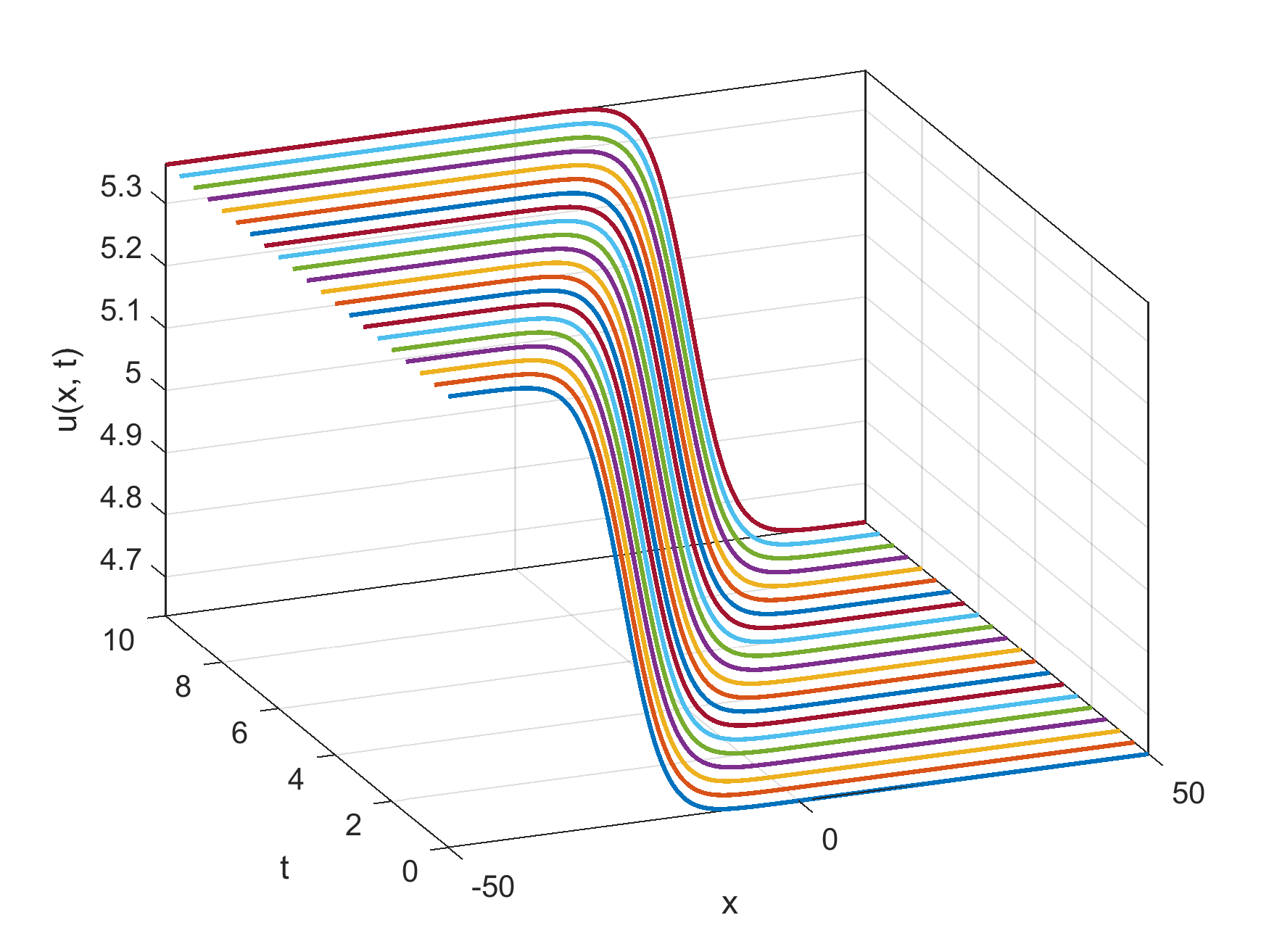}}
\centerline{\text{(a)~Numerical solution}}
\end{minipage}
\begin{minipage}[b]{0.5\linewidth}
\centering
\centerline{\includegraphics[width=\linewidth,height=\textheight,keepaspectratio]{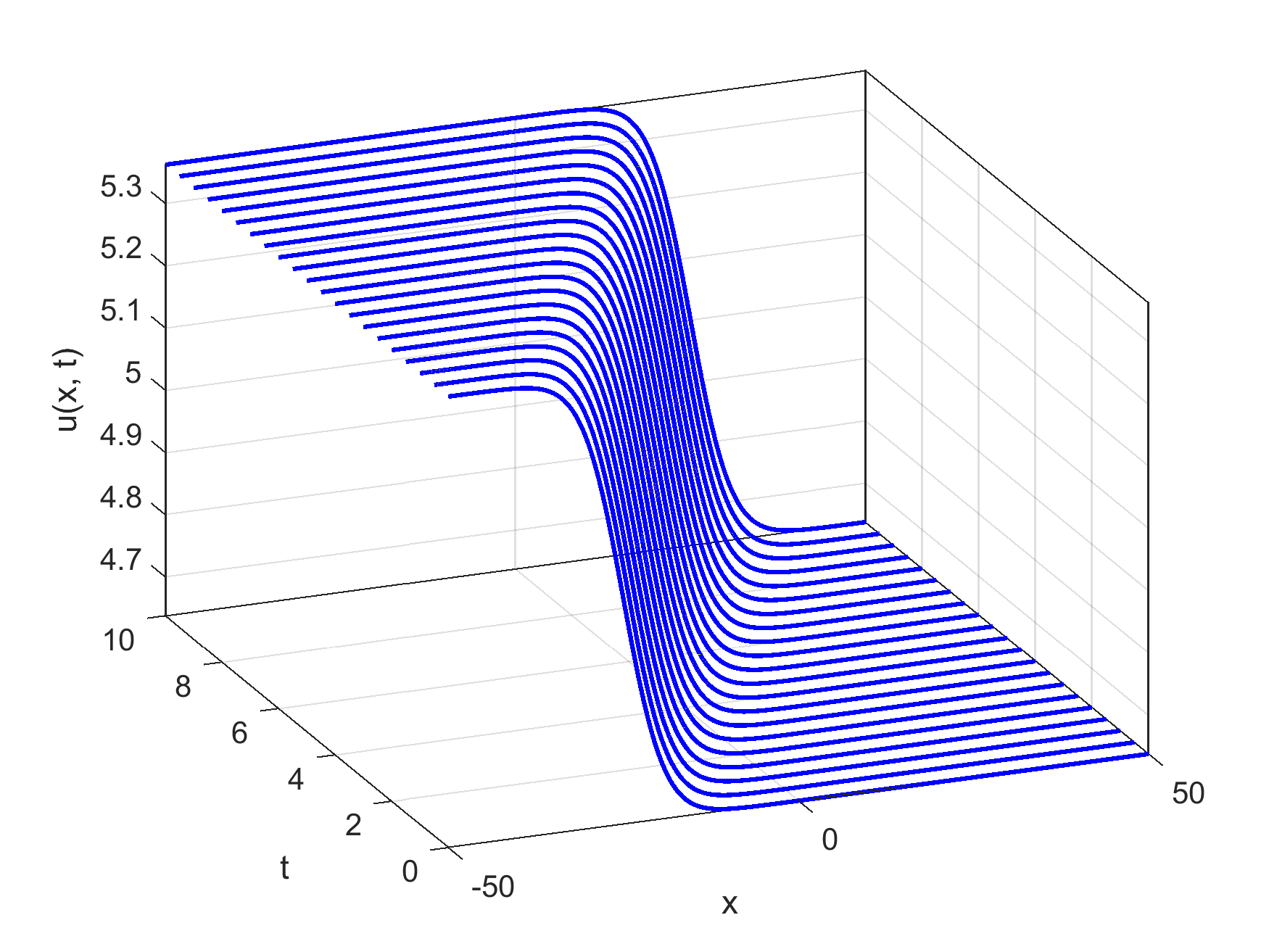}}
\centerline{\text{(b)~Exact solution}}
\end{minipage}
\caption{\footnotesize{Numerical solution vs exact solution for various time $t\in[0,10]$.}}
\label{fig:space and time cgce}
\end{figure}

\begin{figure}[H]
\begin{minipage}[b]{0.5\linewidth}
\centering
\centerline{\includegraphics[width=\linewidth,height=\textheight,keepaspectratio]{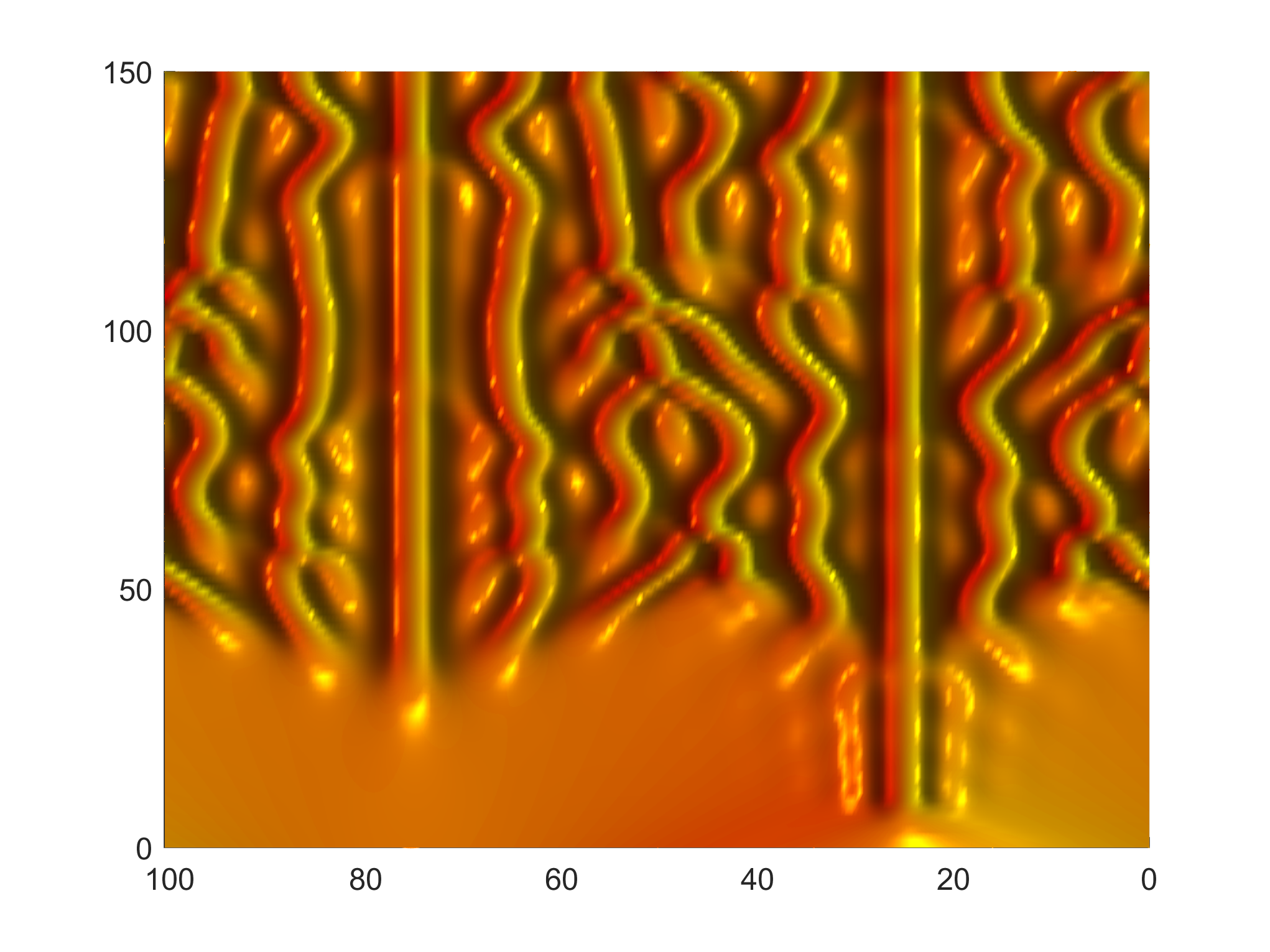}}
\centerline{\text{(a)\ $t\in [0,150], N=256$ and $k=\frac{1}{4}$}}
\end{minipage}
\begin{minipage}[b]{0.5\linewidth}
\centering
\centerline{\includegraphics[width=\linewidth,height=\textheight,keepaspectratio]{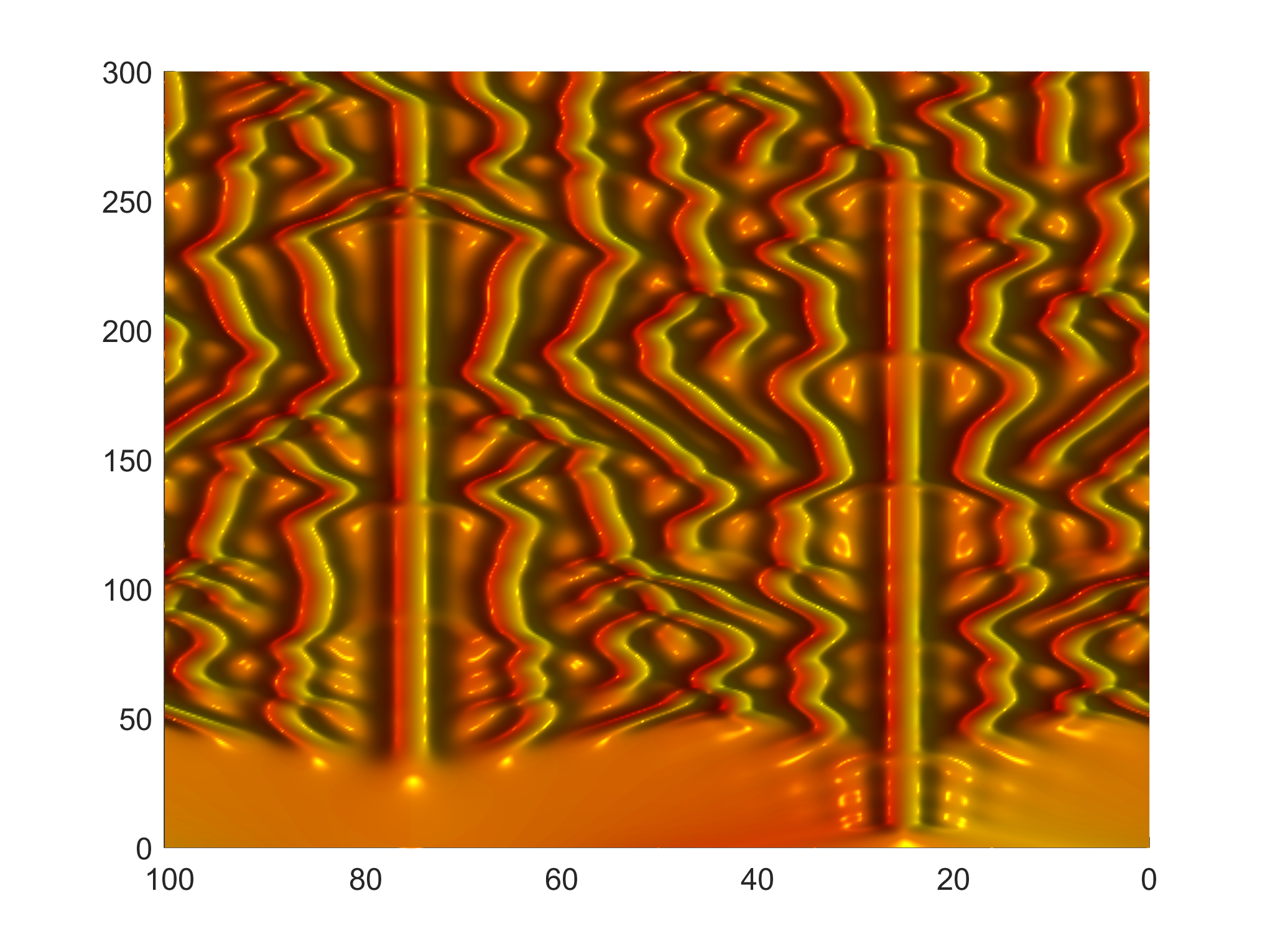}}
\centerline{\text{(b)\ $t\in [0,300], N=512$ and $k=\frac{1}{8}$}}
\end{minipage}
\caption{\footnotesize{Chaotic solution profile of the component $u(x,t)$ at various time obtained via \textup{IMEXRK4}.}}
\label{fig:chaotic}
\end{figure}

\begin{figure}[H]
\begin{minipage}[b]{0.5\linewidth} 
\centering
\centerline{\includegraphics[width=\linewidth,height=\textheight,keepaspectratio]{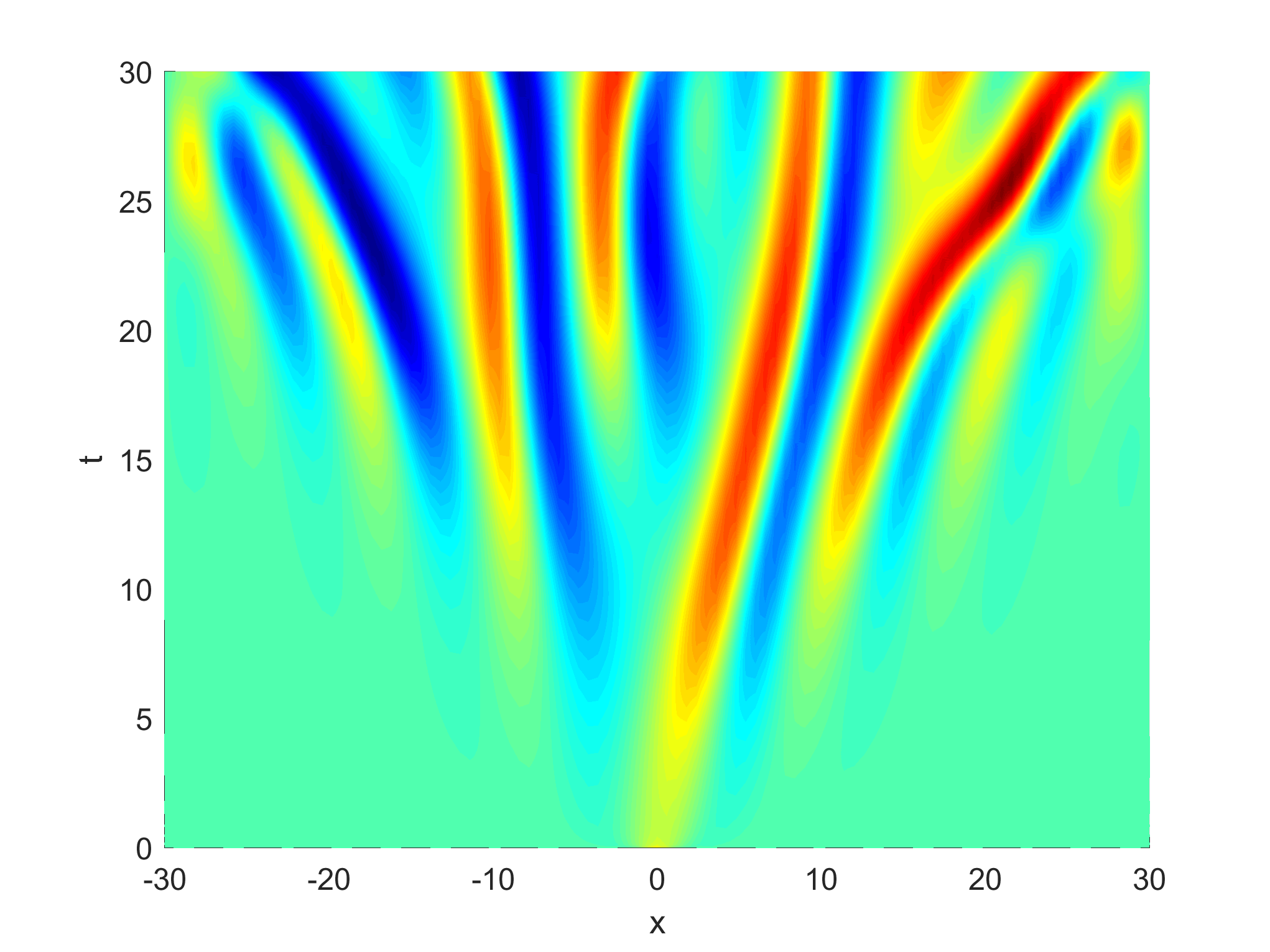}}
\centerline{\text{(a)\ Aerial view}}
\end{minipage}
\begin{minipage}[b]{0.5\linewidth}
\centering
\centerline{\includegraphics[width=\linewidth,height=\textheight,keepaspectratio]{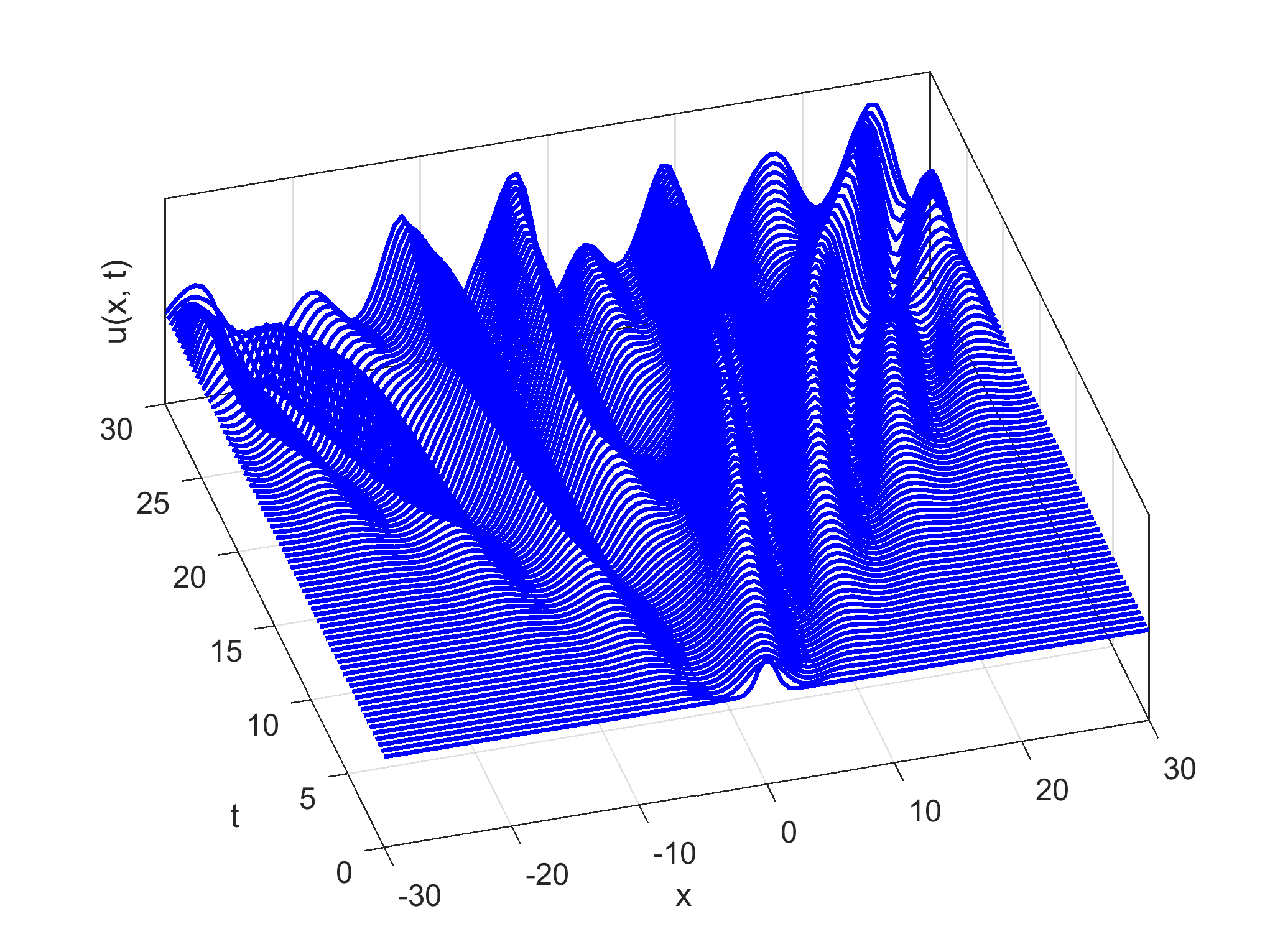}}
\centerline{\text{(b)\ \textup{3D} view}}
\end{minipage}
\caption{\footnotesize{The aerial and \textup{3D} chaotic solution profile of the \textup{KSE} with $N=101,~k=0.1$ at $t\in [0,30]$ for Example 3.}}
\label{fig:pulsev}
\end{figure}

\begin{figure}[H]
\begin{minipage}[b]{0.5\linewidth}
\centering
\centerline{\includegraphics[scale=0.5]{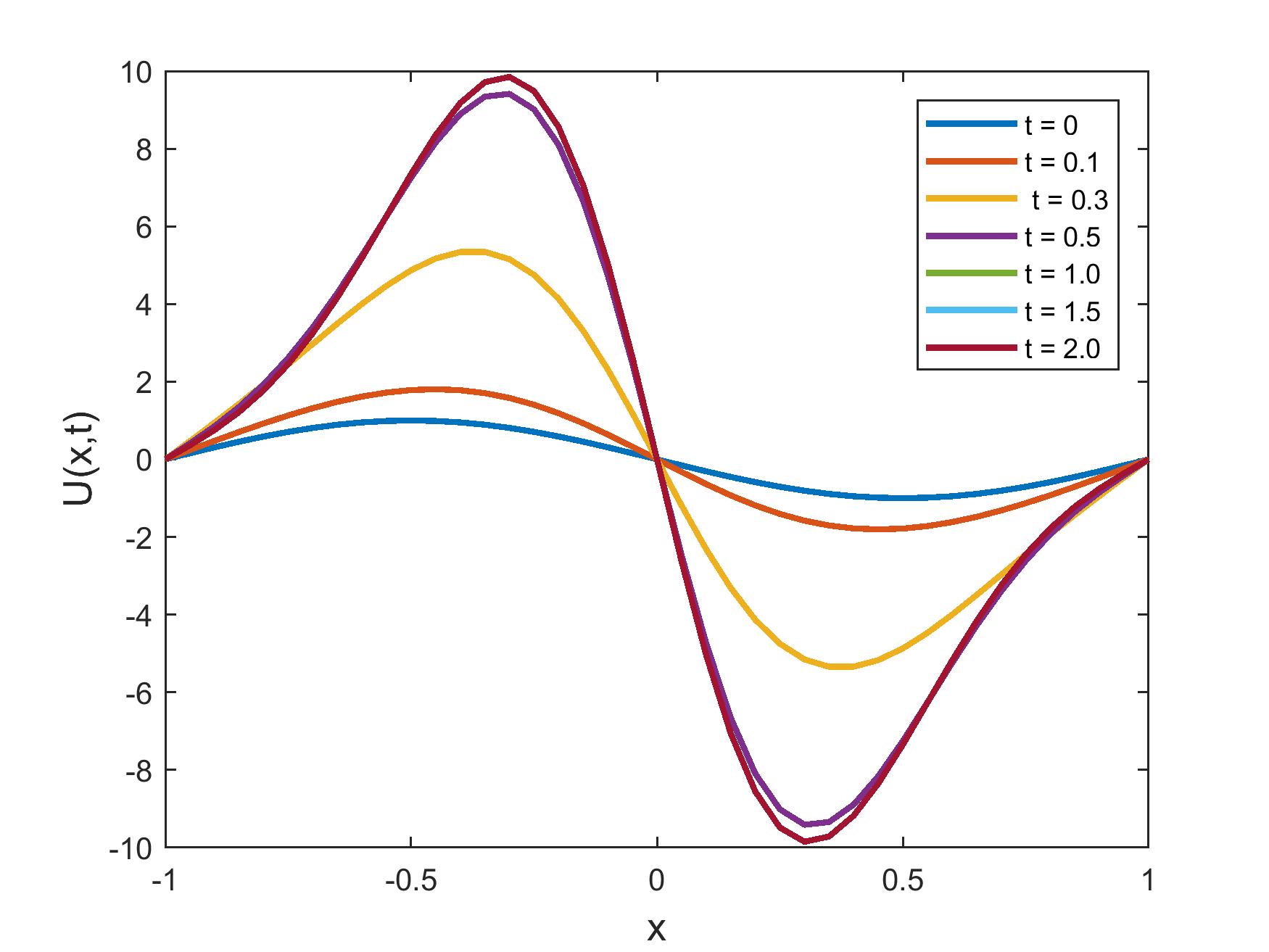}}
\centerline{\text{(a)\ $\beta=0.4/\pi^2$}}
\end{minipage}
\begin{minipage}[b]{0.5\linewidth}
\centering
\centerline{\includegraphics[scale=0.5]{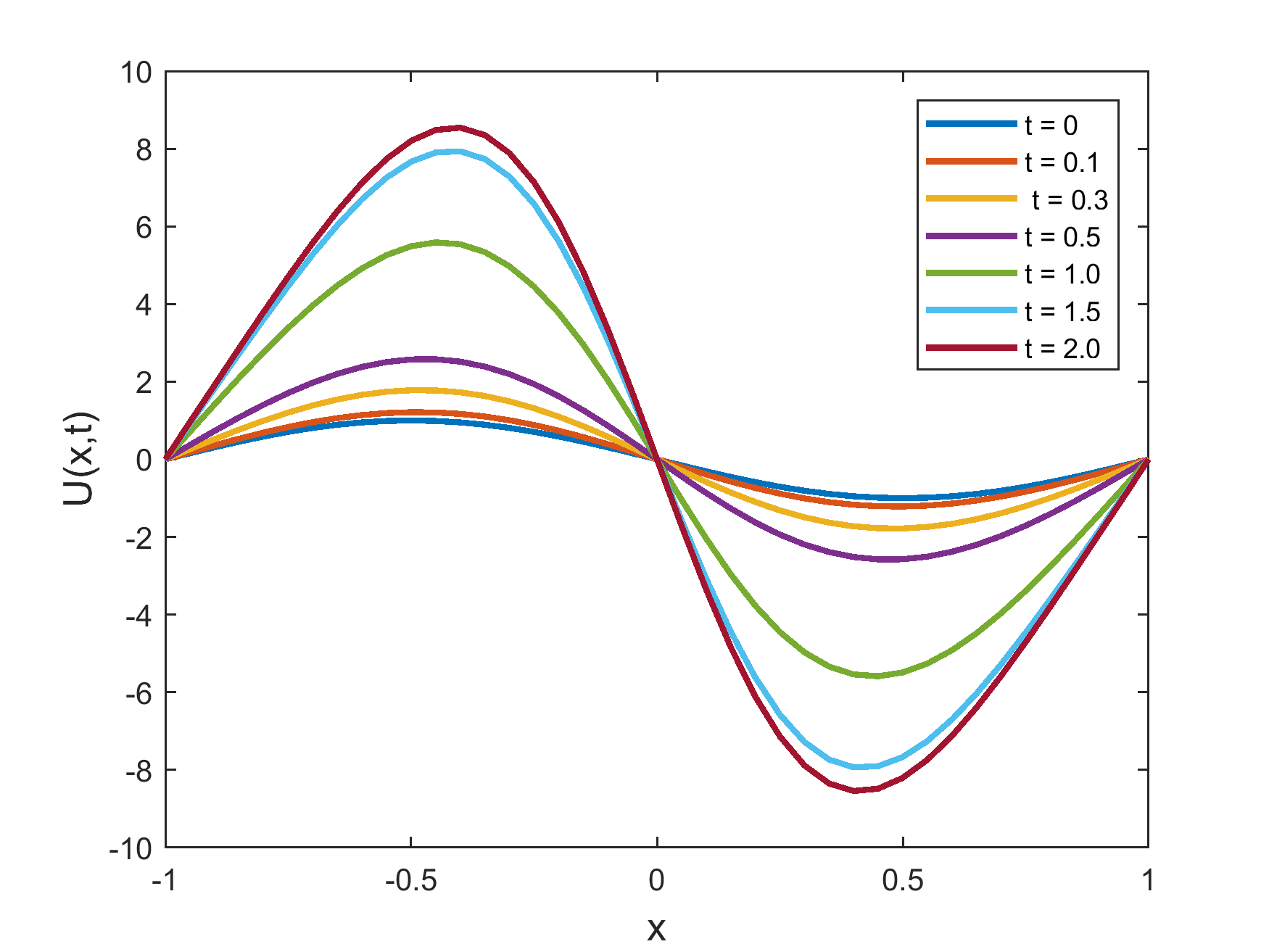}}
\centerline{\text{(b)\  $\beta=0.6/\pi^2$}}
\end{minipage}
\begin{minipage}[b]{0.5\linewidth}
\centering
\centerline{\includegraphics[scale=0.5]{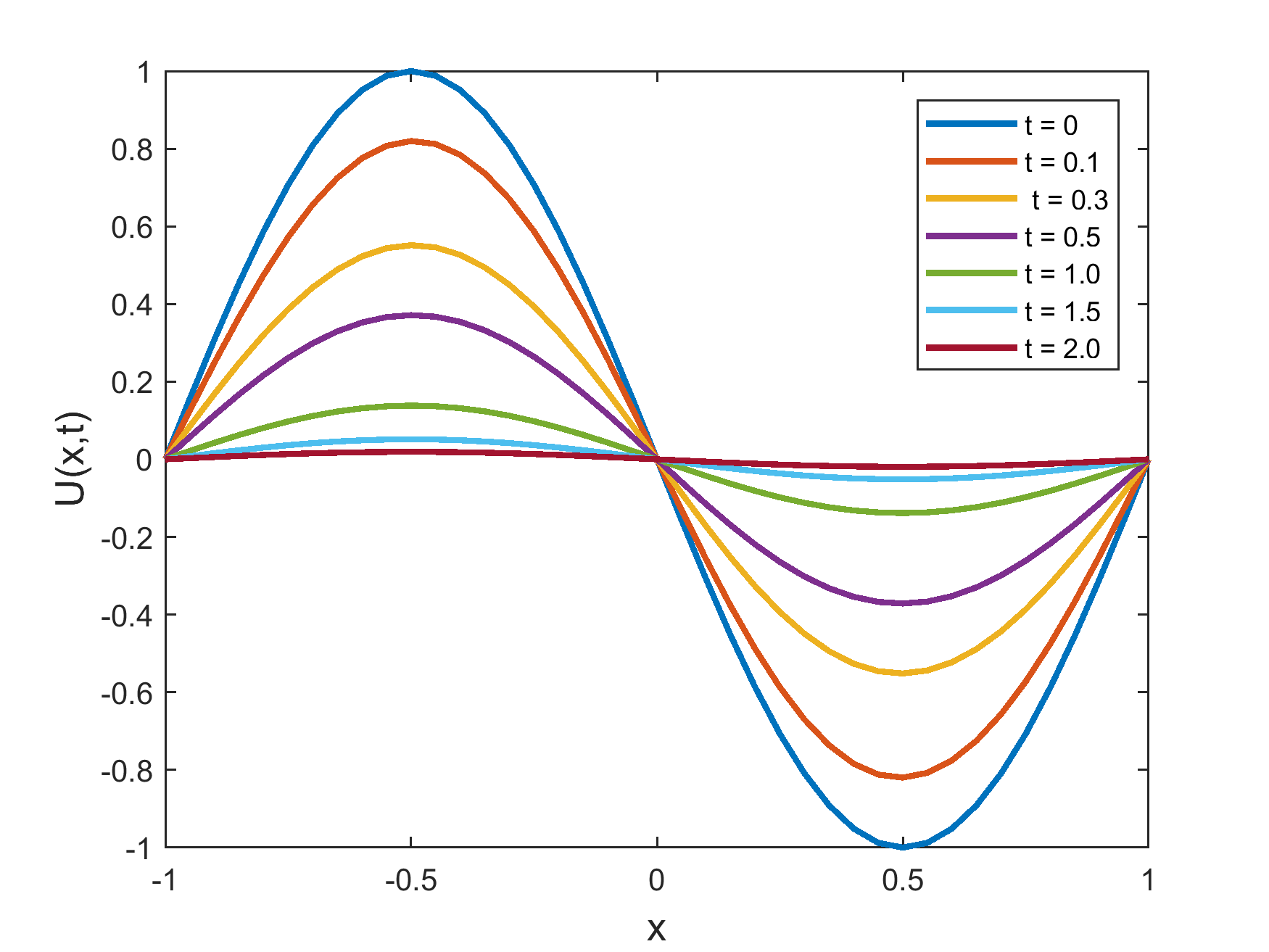}}
\centerline{\text{(c)\  $\beta=0.8/\pi^2$}}
\end{minipage}
\begin{minipage}[b]{0.5\linewidth}
\centering
\centerline{\includegraphics[scale=0.5]{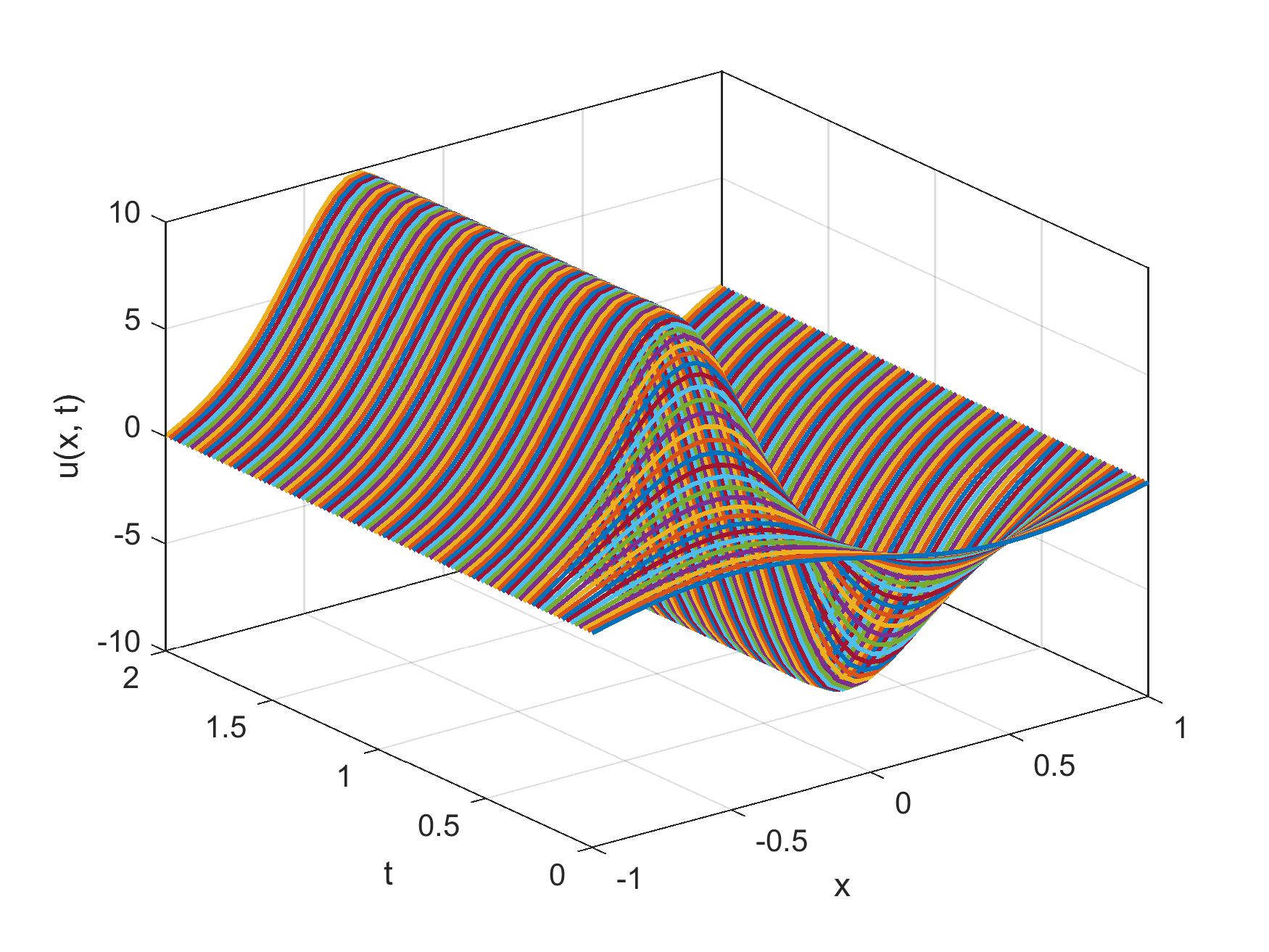}}
\centerline{\text{(d)\  $\beta=0.4/\pi^2$ and $t\in[0,2]$}}
\end{minipage}
\caption{\footnotesize{Space-time evolution profile of the component $u(x,t)$ obtained via \textup{IMEXRK4} scheme with various values of $\beta$ at different time levels.}}
\label{fig:mitt}
\end{figure}

In order to investigate the order of accuracy and computational efficiency of the proposed scheme \textup{IMEXRK4} for solving \textup{KSE}, 
a numerical test on Example 1 was performed. In the computation, a simulation was run up to $t=2.0$ by initially setting $h=4$ and $k=0.025$ then reduced 
both of them by a factor of $2$ in each refinement. The maximum error and rates of convergence are listed in Table \ref{cgce}. As can be seen from Table \ref{cgce} 
that the computed convergence rates of the proposed schemes apparently demonstrate the expected fourth-order accuracy in both time and space.
\begin{table}[H]
\caption{The maximum error, rates of convergence, and CPU time of \textup{IMEXRK4} for Example 1 at $t=2.0$} 
\begin{center} \footnotesize
\resizebox{\linewidth}{!}{
\begin{tabular}{lllllllllllllllllllllllll}
 \hline 
 h \& k&&& 4 \& 0.025 &&&2 \& 0.0125&&&1 \& 0.00625&&&0.5 \& 0.003125\\ 
\cline{1-13}
$\norm{\cdot}_\infty$&&& 6.157E-03&&& 3.775E-04&&& 2.396E-05&&& 1.461E-06\\
Order&&& \quad - &&&4.0278&&& 3.9777&&& 4.0359\\ 
CPU(s)&&&0.2004&&&0.5011&&& 1.6728&&&8.0996\\ 
\hline
\end{tabular}}
\end{center} 
\label{cgce} 
\end{table}
In order to visualize the space and time rates of convergence of the proposed scheme, we illustrated them in Fig.\,\ref{fig:spacetime cgce} with log-log scale graph. 
From Fig.\,\ref{fig:spacetime cgce}, it can be seen that the slopes of the regression line for maximum errors both in time and spatial directions are 
close to four, which corresponds to the fourth-order scheme both in time and space.
We ran another sets of experiment in Example 1 with $h=0.5, k=0.01$ until time $t=10$ and captured the \textup{3D} view of the solution profile in 
Fig.\,\ref{fig:space and time cgce} (a). From Fig.\,\ref{fig:space and time cgce} (a)-(b), we can see the the solution obtained via using the proposed scheme 
is close to the exact solution.
In Table \ref{tablegre}, we compare the scheme \textup{IMEXRK4} with earlier schemes: SBSC\cite{ZMP}, QBSC\cite{RGA}, and LBM\cite{LHM} by listing GRE of 
$u(x,t)$ with $N=200,~k=0.01$ at different time levels $t\leq 12.$ Table \ref{tablegre} clearly indicates that the numerical results obtained via proposed scheme 
are more accurate than that obtained via SBSC\cite{ZMP}, QBSC\cite{RGA} and LBM\cite{LHM}.
\begin{table}[H]
\caption{Comparison of GRE at different time with $N=200$ and $k=0.01$ for Example 1.} 
\begin{center} \footnotesize
\resizebox{\linewidth}{!}{
\begin{tabular}{llllllllllllllllll}
\hline
&&&Time$(t)$ &&& 6 &&& 8 &&& 10 &&& 12\\
\cline{4-16}
IMEXRK4&&&GRE&&& 7.624E-08&&& 8.092E-08&&& 8.589E-08&&& 3.188E-07\\\\
\cline{4-16}
&&&Time$(t)$ &&& 6 &&& 8 &&& 10 &&& 12\\
\cline{4-16}
SBSC &&&GRE&&& 1.625E-07&&& 1.940E-07&&& 2.229E-07&&& 5.314E-07\\\\
\cline{4-16}
&&&Time$(t)$ &&& 6 &&& 8 &&& 10 &&& 12\\
\cline{4-16}
QBSC&&&GRE&&& 6.509E-06&&& 7.132E-06&&& 7.310E-06&&& 8.776E-06\\\\
\cline{4-16}
&&&Time$(t)$ &&& 6 &&& 8 &&& 10 &&& 12\\
\cline{4-16}
LBM&&&GRE&&& 7.881E-06&&& 9.532E-06&&& 1.089E-05&&& 1.179E-05\\\\
\hline
\end{tabular}}
\end{center} 
\label{tablegre} 
\end{table}
\noindent {\bf Example 2 (\textup{KSE} with periodic boundary conditions):}\ In this example, we consider the \textup{KSE} with $\alpha=1$ and $\beta=1$ over a domain 
$\Omega=[0,32\pi]$ along with periodic boundary conditions and following initial condition:
\begin{flalign}
&u(x,0)=\cos(\frac{x}{16})\big(1+\sin(\frac{x}{16})\big).&
\end{flalign}\label{eq:prob2}
Two sets of numerical experiments on Example 2 were conducted. In the first sets of the experiment, the order of accuracy in the temporal direction of the 
proposed method with periodic boundary conditions was checked by running an experiment until time $T = 10.0$ with $N=256$. Initially, $k=\frac{1}{2}$ was set and 
repeatedly halved it at each time and numerical results are captured in Table \ref{cgcekas}. The error values $E_k$ listed for \textup{IMEXRK4}
scheme in Table \ref{cgcekas} are again calculated through a maximal difference between each simulation. From Table \ref{cgcekas}, 
it can be seen that the proposed scheme is able to achieve the expected fourth-order accuracy in time with periodic boundary conditions.
\begin{table}[H]
\caption{The maximum error, time rates of convergence, and CPU time of \textup{IMEXRK4} for Example 2 with $N=256$ at $t=10.0$} 
\begin{center} \footnotesize
\resizebox{\linewidth}{!}{
\begin{tabular}{lllllllllllllllllllllllll}
 \hline
 k&&& 1/4 &&&1/8&&&1/16&&&1/32\\ 
\cline{1-13}
$E_k$&&& 9.031E-04&&& 6.291E-05&&& 3.922E-06&&& 2.442E-07\\
Order&&& \quad - &&&3.8436&&& 4.0034&&& 4.0052\\ 
CPU(s)&&&1.1791&&&2.1269&&& 3.3537&&&5.2132\\ 
\hline
\end{tabular}}
\end{center} 
\label{cgcekas} 
\end{table}
To better understand the applicability of the \textup{IMEXRK4} scheme to simulate the long time behavior of the \textup{KSE}, 
a second sets of experiment was conducted on Example 2 until the long time $t=150$ with $N=256,~~k=\frac{1}{4}$ and $t=300$ with $N=512,~~k=\frac{1}{8}$. 
The chaotic solution profiles of the component $u(x, t)$ for $t\in[0,150]$ and for $t\in[0,300]$ were captured in Fig.\,\ref{fig:chaotic}\,(a) and (b) respectively. 
Chaotic solution profile corresponds for $t\leq 150$ in the Fig.\,\ref{fig:chaotic}\,(a) are in good agreement with the results depicted in \cite{KTN} -- therefore, we are 
confident that the profile corresponds for $t\leq 300$ is correct and reliable.
\noindent {\bf Example 3 (\textup{KSE} with Gaussian initial condition):}\ Here, we consider the \textup{KSE} with $\alpha=1$ and $\beta=1$ which exhibits chaotic 
behavior over a finite domain $\Omega=[-30,30]$ with homogeneous Dirichlet boundary conditions and the following Gaussian initial condition
\begin{flalign}
&u(x,0)=\exp(-x^2).&\label{eq:gas}
\end{flalign}
In what follows, again on Example 3 two sets of numerical experiments were performed. In the first sets of the experiment, the order of accuracy in the temporal 
direction of the proposed scheme with homogeneous Dirichlet boundary conditions was examined by running an experiment until time T = 1.0 with fixed $N=101$. 
Again the error values reported in Table \ref{cgcegas} are calculated through a maximal difference between each simulation which is obtained by repeatedly halving an 
initial time step size $k=\frac{1}{2}$ at each time. From the results, it can be seen that the proposed scheme is able to achieve the expected fourth-order accuracy in 
time.
\begin{table}[H]
\caption{The maximum error, time rates of convergence, and CPU time of \textup{IMEXRK4} for Example 3 with $N=101$ at $t=1.$} 
\begin{center} \footnotesize
\resizebox{\linewidth}{!}{
\begin{tabular}{lllllllllllllllllllllllll}
 \hline
 k&&& 0.01/2 &&&0.01/4&&&0.01/8&&&0.01/16\\ 
\cline{1-13}
$E_k$&&& 2.723E-08&&& 1.976E-09&&& 1.324E-10&&& 8.613E-12\\
Order&&& \quad - &&&3.7847&&& 3.8995&&& 3.9422\\ 
CPU(s)&&&0.5543&&&1.0044&&& 1.9298&&&3.4635\\ 
\hline
\end{tabular}}
\end{center} 
\label{cgcegas} 
\end{table}
In the second sets of experiment, the chaotic behavior of the component $u(x,t)$ is simulated for Gaussian initial condition. 
The simulations are accomplished in $t\in[0, 30]$ with the parameters $N= 101$ and $k=0.1$ and captured in Fig. \ref{fig:pulsev} with aerial and \textup{3D} views. 
From Fig. \ref{fig:pulsev}, we can see clearly that the result shows same behavior as reported in \cite{LMD,RGA}.
\noindent {\bf Example 4 (\textup{KSE} with different values of $\beta$):}\ The \textup{KSE} with $\alpha=1$ and different values of $\beta$ over a finite domain 
$\Omega=[-1,1]$ along with homogeneous Dirichlet boundary conditions and the following initial condition is considered here:
\begin{flalign}
&u(x,0)=-\sin(\pi x).&\label{eq:gray}
\end{flalign}
For the empirical convergence analysis in the temporal direction of the proposed scheme,  we ran an experiment until time T = 1.0 with $\beta=1.1$ and fixed $h=0.05$. 
The error values reported in Table \ref{cgcemitt} are computed through a maximal difference between each simulation which is obtained by repeatedly halving an initial 
time step size $k=0.005$ at each time. From the results, it can be seen that the proposed scheme is able to exhibit the expected fourth-order accuracy in time. 
\begin{table}[H]
\caption{The maximum error, time rates of convergence, and CPU time of \textup{IMEXRK4} for Example 4 with $\beta=1.1,~h=0.05$ at $t=1.$} 
\begin{center} \footnotesize
\resizebox{\linewidth}{!}{
\begin{tabular}{lllllllllllllllllllllllll}
 \hline 
 k&&& 0.005/2 &&&0.005/4&&&0.005/8&&&0.005/16\\ 
\cline{1-13}
$E_k$&&& 1.431E-08&&& 9.7926E-10&&& 6.532E-11&&& 3.320E-12\\
Order&&& \quad - &&&3.8692&&& 3.9060&&& 4.2983\\ 
CPU(s)&&&1.0360&&&2.1329&&& 3.7341&&&5.6047\\ 
\hline
\end{tabular}}
\end{center} 
\label{cgcemitt} 
\end{table}
The space-time evolution profile of $u(x,t)$ via \textup{IMEXRK4} scheme with different values of $\beta$, $h=0.05,~k=0.001$ at different time levels were 
depicted in Fig.\,\ref{fig:mitt}. The results in Fig.\,\ref{fig:mitt} clearly exhibit good agreement with results reported in \cite{RGA, AVM}.
\section{Conclusions}

This manuscript introduced a fourth-order scheme both in time and space to solve the KSE. The proposed method utilized a compact fourth-order finite difference scheme for a spatial discretization and the fourth-order Runge-Kutta based implicit-explicit scheme for time discretization. A Compact finite difference scheme is used to transform the KSE to a system of ordinary differential equations (ODEs) in time, and then, fourth-order time integrator is implemented to solve the resulting ODEs. Calculation of local truncation error and an empirical convergence analysis exhibited the fourth-order accuracy of the proposed scheme. The performance and applicability of the scheme have been investigated by testing it on several test problems. The computed numerical solutions maintain a good accuracy compared with the exact solution. In addition, the numerical results exhibited that the proposed scheme provides better accuracy in comparison with other existing schemes.

\bibliographystyle{unsrt}  


\end{document}